%% file: Main_arXiv.tex
\documentclass[a4paper,
               11pt,
               pdftex,
               normalheadings,
               headsepline,
               footsepline,
               onecolumn,
               headinclude,
               footinclude,
               DIV14,
               abstracton]
{scrartcl}

\usepackage
{
    graphicx,
    amssymb,
    amsmath,
    amsthm,
    xcolor,
    dsfont,
    algpseudocode,
    authblk,
}

\usepackage{blindtext} 

\usepackage[bf]{caption}
\usepackage[colorlinks,
            pdffitwindow=false,
            plainpages=false,
            pdfpagelabels=true,
            pdfpagemode=UseOutlines,
            pdfpagelayout=SinglePage,
            bookmarks=false,
            colorlinks=true,
            hyperfootnotes=false,
            linkcolor=blue,
            citecolor=green!50!black]
{hyperref}

\usepackage{enumerate}
\usepackage{algorithm}

\DeclareMathAlphabet{\mathpzc}{OT1}{pzc}{m}{it}

\newtheorem{theorem}{Theorem}[section]

\newtheorem{definition}[theorem]{Definition}
\newtheorem{remark}[theorem]{Remark}

\makeatletter
\renewcommand*\env@matrix[1][*\c@MaxMatrixCols c]{%
  \hskip -\arraycolsep
  \let\@ifnextchar\new@ifnextchar
  \array{#1}}
\makeatother

\begin{document}

\title{Multiobjective Optimal Control Methods for Fluid Flow Using Reduced Order Modeling}
\author[*]{Sebastian Peitz}
\author[**]{Sina Ober-Bl{\"o}baum}
\author[*]{Michael Dellnitz}
\affil[*]{\normalsize Department of Mathematics, University of Paderborn, Warburger Str.~100, D-33098 Paderborn}
\affil[**]{Department of Engineering Science, University of Oxford, Parks Road, Oxford OX1 3PJ, UK}

\maketitle

\begin{abstract}
In a wide range of applications it is desirable to optimally control a dynamical system with respect to concurrent, potentially competing goals. This gives rise to a multiobjective optimal control problem where, instead of computing a single optimal solution, the set of optimal compromises, the so-called Pareto set, has to be approximated. When the problem under consideration is described by a partial differential equation (PDE), as is the case for fluid flow, the computational cost rapidly increases and makes its direct treatment infeasible. Reduced order modeling is a very popular method to reduce the computational cost, in particular in a multi query context such as uncertainty quantification, parameter estimation or optimization. In this article, we show how to combine reduced order modeling and multiobjective optimal control techniques in order to efficiently solve multiobjective optimal control problems constrained by PDEs. We consider a global, derivative free optimization method as well as a local, gradient based approach for which the optimality system is derived in two different ways. The methods are compared with regard to the solution quality as well as the computational effort and they are illustrated using the example of the two-dimensional incompressible flow around a cylinder. 
\end{abstract}

\pagestyle{myheadings}
\thispagestyle{plain}

\section{Introduction}
\label{sec:Introduction}
\input{01_Introduction.tex}

\addtocounter{equation}{1}
\section{Problem formulation}
\label{sec:Problem}
\input{02_Problem.tex}

\section{Multiobjective Optimal Control}
\label{sec:MOC}
\input{03_MOCP.tex}

\section{Reduced Order Modeling}
\label{sec:ROM}
\input{04_ROM.tex}

\section{Results}
\label{sec:Results}
\input{05_Results.tex}

\section{Conclusion}
\label{sec:Conclusion}
\input{06_Conclusion.tex}

\textbf{Acknowledgement:} This research was funded by the Leading-Edge Cluster ''Intelligent Technical Systems OstWestfalenLippe'' (its OWL) and the DFG Priority Programme 1962 ''Non-smooth and Complementarity-based Distributed Parameter Systems''. Calculations leading to the results presented here were performed on resources provided by the Paderborn Center for Parallel Computing ($PC^2$).

\bibliographystyle{alpha}
\bibliography{Bib_update}

\input{Appendix.tex}

\end{document}

%% file: 01_Introduction.tex
In many applications from industry and economy, one is interested in simultaneously optimizing several criteria. For example, in transportation one wants to reach a destination as fast as possible while minimizing the energy consumption. This example illustrates that in general, the different objectives contradict each other. Therefore, the task of computing the set of optimal compromises between the conflicting objectives, the so-called \emph{Pareto set}, arises. This leads to a multiobjective optimization problem (MOP) or, if the control variable is a function, a multiobjective optimal control problem (MOCP). Based on the knowledge of the Pareto set, a \emph{decision maker} can use this information either for improved system design or for changing control parameters during operation as a reaction on external influences or changes in the system state itself.

Multiobjective optimization is an active area of research. Different approaches exist to address MOPs, e.g.~deterministic approaches \cite{Mie12, Ehr05}, where ideas from scalar optimization theory are extended to the multiobjective situation. In many cases, the resulting solution method involves solving multiple scalar optimization problems consecutively. Continuation methods make use of the fact that under certain smoothness assumptions the Pareto set is a manifold that can be approximated by continuation methods known from dynamical systems theory \cite{Hil01}. Another prominent approach is based on evolutionary algorithms \cite{AJG05, CCLvV07}, 
where the underlying idea is to evolve an entire set of solutions (population) during the optimization process. Set oriented methods provide an alternative deterministic approach to the solution of MOPs. Utilizing subdivision techniques, the desired Pareto set is approximated by a nested sequence of increasingly refined box coverings \cite{DSH05, SWO+13}.

When dealing with control functions, a multiobjective optimal control problem (MOCP) needs to be solved. Similar to MOPs, ideas from scalar optimal control theory (see e.g.~\cite{HPU+08} for an overview of optimal control methods for PDE-constrained problems) can be extended to take into account multiple objectives. By applying a direct method, the MOCP is transformed into a high-dimensional, nonlinear MOP such that the methods mentioned before can be applied. 
Another approach is based on the transformation of the MOCP into a sequence of scalar optimal control problems and the use of well established optimal control techniques for their solution. Examples for MOCPs can be found e.g. in \cite{DEF+16, LHDI10, ORF12, SWO+13} with ODE constraints. 
In \cite{LKBM05} as well as \cite{ARFL09}, nonlinear PDE constraints are taken into account but the model is treated as a black box, i.e.~no special treatment of the constraints is required. The first articles explicitly taking into account PDE constraints are \cite{IUV13, ITV15}, where multiobjective optimal control problems are solved with a weighted sum approach and model order reduction techniques subject to linear and semilinear PDE constraints, respectively. Fluid flow applications have been considered in \cite{OBPG15, PD15}.

All approaches to MOP / MOCP have in common that a large number of function evaluations is typically required. Thus, the direct computation of the Pareto set can quickly become numerically infeasible. This is frequently the case for problems described by (nonlinear) partial differential equations such as the Navier-Stokes equations. Standard optimization methods for PDEs \cite{HPU+08} often make use of a discretization by finite elements, finite volumes or finite differences which results in a high-dimensional system of ODEs. In a multi query context (such as optimization, parameter identification, etc.), this approach often exceeds the limits of today's computing power. Hence, one aims for methods which reduce the computational cost significantly. This can be achieved by approximating the PDE by a reduced order model of low dimension.

In recent years, major achievements were made in the field of reduced order modeling \cite{SVR05, BMS05}. In fact, different methods for creating low dimensional models exist for linear systems (e.g.~\cite{WP02}) as well as for nonlinear systems (e.g.~\cite{BMNP04, GMN+07}), see \cite{ASG01} for a survey. Many researchers focus their attention to the Navier-Stokes equations (e.g.~\cite{CBS05, Row05, XFB+14}), where \emph{Proper Orthogonal Decomposition} (POD) \cite{HLB98} has proven to be a powerful tool. In the context of fluids, this technique has first been introduced by Lumley \cite{Lum67} to identify coherent structures. Reduced order models using POD modes and the \emph{method of snapshots} go back to Sirovich \cite{Sir87}.

Due to the wide spectrum of potential applications (optimal mixing, drag reduction, HVAC (\emph{Heating, Ventilation and Air Conditioning}), etc.) and the progress in computational capabilities, a lot of research is devoted to the control of the Navier-Stokes equations, either directly \cite{GeH89, Lac14} or via reduced order modeling \cite{GPT99, Fah00, IR01, BCB05, BC08, PD15}. In many cases, the energy consumption plays an important role. Thus, ideally, one wants to consider the control cost in addition to the main objective in order to choose a good trade-off between the achievement of the desired objective and the respective control cost. Since this causes a considerable computational effort in case of systems described by the Navier-Stokes equations, we show in this article how reduced order modeling and multiobjective optimal control methods can be combined to solve MOCPs with nonlinear PDE constraints. Using model order reduction via POD and Galerkin projection, we present several methods to compute the Pareto set for the conflicting objectives \emph{minimization of flow field fluctuations} and \emph{control cost} for the two-dimensional flow around a cylinder at $Re = 200$. We discuss the advantages and disadvantages of the different approaches and comment on the respective computational cost. In fact, it is shown that the different methods strongly differ in their respective efficiency and that gradient based approaches have a better performance, provided that the gradient is sufficiently accurate. We also show that the choice of the algorithm and the corresponding particular numerical realization of the control function lead to different optimal control behavior.

The article is organized as follows. In Section~\ref{sec:Problem} we present the problem setting and formulate the MOCP. After an introduction to multiobjective optimal control in Section~\ref{sec:MOC}, the reduced order model and the resulting reduced MOCP are derived in Section~\ref{sec:ROM}. We then present our results in Section~\ref{sec:Results} and draw a conclusion in Section~\ref{sec:Conclusion}.

%% file: 02_Problem.tex
\begin{figure}
	\centering
	\includegraphics[width = 0.65\textwidth]{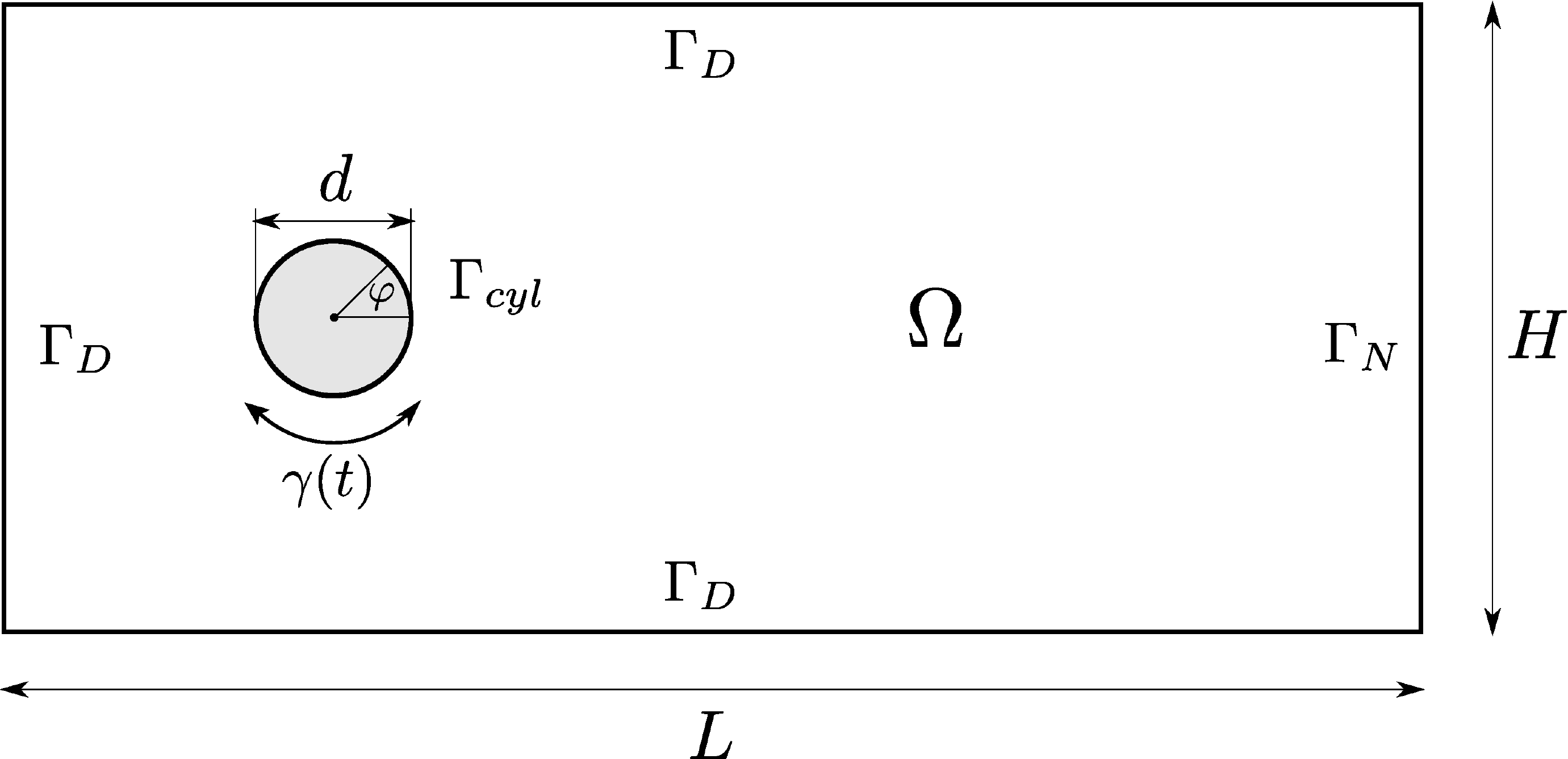}
	\caption{Sketch of the domain $\Omega \subset \mathbb{R}^2$. The length is $L = 25d$, the height $H = 15d$, and the cylinder center is placed at $(5d, 7.5d)$.}
	\label{fig:Domain}
\end{figure}
The two-dimensional, viscous flow around a cylinder is one of the most extensively studied cases for separated flows in general as well as for flow control problems \cite{GPT99, GBZI04, BCB05}. In this paper, we consider the laminar case described by the incompressible Navier-Stokes equations at a Reynolds Number $Re = \frac{\boldsymbol{U}_{\infty} d}{\nu} = 200$ computed with respect to the far field velocity $\boldsymbol{U}_{\infty}$ (throughout the paper, we will use bold symbols for vector valued quantities), the kinematic viscosity $\nu$ and the cylinder diameter $d$:
\begin{align}
	\frac{\partial \boldsymbol{U}(\boldsymbol{x},t)}{\partial t} + \left(\boldsymbol{U}(\boldsymbol{x},t) \cdot \nabla \right) \boldsymbol{U}(\boldsymbol{x},t) &= -\frac{\nabla p(\boldsymbol{x},t)}{\rho} + \frac{1}{Re} \nabla^2 \boldsymbol{U}(\boldsymbol{x},t), \tag{2.1a} \label{eq:NSE1} \\
	\nabla \cdot \boldsymbol{U}(\boldsymbol{x},t) &= 0,  \tag{2.1b} \label{eq:NSE2} \\
	\left( \boldsymbol{U}(\boldsymbol{x}, 0), p(\boldsymbol{x}, 0) \right) &= \left(\boldsymbol{U}_0(\boldsymbol{x}), p_0(\boldsymbol{x})\right), \tag{2.1c} \label{eq:NSE_IC} \\
	\text{for}~\boldsymbol{x} \in \Omega,\ &t \in [t_0, t_e], \notag
\end{align}
where $\boldsymbol{U} \in H^2(\Omega \times [t_0,t_e], \mathbb{R}^2)$ is the two-dimensional fluid velocity and $p \in H^1(\Omega \times [t_0,t_e], \mathbb{R})$ the pressure. 
$H^k$ is the standard Sobolev space $W^{k,2}$ (cf. e.g.~\cite{HPU+08}). 
The domain (cf.~Figure~\ref{fig:Domain}) is denoted by $\Omega$. 
We impose Dirichlet boundary conditions at the inflow as well as the upper and lower walls $\Gamma_D$. At the outflow $\Gamma_N$, we impose a standard \emph{no shear stress} condition \cite{HRT96}:
\begin{align}
	\boldsymbol{U}(\boldsymbol{x},t) &= (\boldsymbol{U}_{\infty}, 0) &\text{for}~\boldsymbol{x} \in \Gamma_{D}, \tag{2.1d} \label{eq:BCD} \\
	p(\boldsymbol{x},t)\ \boldsymbol{n} &= \frac{1}{Re} \frac{\partial \boldsymbol{U}(\boldsymbol{x},t)}{\partial \boldsymbol{n}} &\text{for}~\boldsymbol{x} \in \Gamma_{N}, \tag{2.1e} \label{eq:BCN}
\end{align}
where $\boldsymbol{n} \in \mathbb{R}^2$ is the outward normal vector of the boundary. On the cylinder $\Gamma_{cyl}$, we prescribe a time-dependent Dirichlet BC such that it performs a rotation around its center with the angular velocity $\gamma(t)$:
\begin{align}
	\boldsymbol{U}(\boldsymbol{x},t) = \frac{d}{2} \gamma(t) \left( \begin{array}{c} -\sin(\varphi) \\ \cos(\varphi) \end{array} \right) \hspace{3cm}\text{ for }\boldsymbol{x} \in \Gamma_{cyl}, \tag{2.1f} \label{eq:BCC}
\end{align}
with $\varphi$ according to Figure \ref{fig:Domain}. The cylinder rotation $\gamma(t) \in L^2([t_0,t_e],\mathbb{R})$ serves as the control mechanism for the flow. The Hilbert space $L^2$ is equipped with the inner product $\left(\boldsymbol{u}, \boldsymbol{v}\right)_{L^2} = \int_{t_0}^{t_e} \boldsymbol{u}(t) \cdot \boldsymbol{v}(t) \, dt$ and the norm $\|\boldsymbol{u}\|_{L^2} = \left(\int_{t_0}^{t_e} \boldsymbol{u}(t) \cdot \boldsymbol{u}(t)\, dt \right)^{1/2}$.

Following \cite{Fah00, BCB05}, we introduce the weak formulation of \eqref{eq:NSE1}. Consider the divergence free Hilbert space of test functions $V = \left\lbrace \boldsymbol{\psi} \in H^1(\Omega \times [t_0, t_e], \mathbb{R}^2) \ | \ \nabla \cdot \boldsymbol{\psi} = 0 \right\rbrace$. Then, a function $\boldsymbol{U} \in H^1(\Omega \times [t_0, t_e], \mathbb{R}^2)$ which satisfies
\begin{align}
	\left(\frac{\partial \boldsymbol{U}}{\partial t} + \left(\boldsymbol{U} \cdot \nabla \right) \boldsymbol{U}, \boldsymbol{\psi} \right) = \left( \frac{p}{\rho}, \nabla \cdot \boldsymbol{\psi} \right) - \left[\frac{p}{\rho} \boldsymbol{\psi} \right] - \frac{1}{Re} \left( \nabla \boldsymbol{\psi}, (\nabla \boldsymbol{U})^\top \right) + \frac{1}{Re} \left[ (\nabla \boldsymbol{U})^\top \boldsymbol{\psi} \right] \tag{2.2} \label{eq:NSE_weak}
\end{align}
for all $\boldsymbol{\psi} \in V$ is called a \emph{weak solution} of \eqref{eq:NSE1}. 
Here, $\left[ \cdot \right]$ is the boundary integral (e.g.~$\left[ \boldsymbol{U} \right] = \int_{\Gamma} \boldsymbol{U}(\boldsymbol{x}) \cdot \boldsymbol{n}\, d\boldsymbol{x}$) and $\left( \cdot, \cdot \right)$ is the inner product for vector-valued quantities (e.g.~$\left( \nabla \boldsymbol{\psi}, (\nabla \boldsymbol{U})^\top \right) = \int_{\Omega} \sum_{i,j} \partial \boldsymbol{\psi}_{i} / \partial x_j \cdot \partial \boldsymbol{U}_{i} / \partial x_j \,d\boldsymbol{x}$). Note that \eqref{eq:NSE2} is automatically satisfied by design of the test function space $V$.

\subsection{Finite volume discretization}
\label{subsec:Problem_solution_method}
\begin{figure}
	\centering
	\parbox[b]{0.49\textwidth}{\centering \includegraphics[width=0.49\textwidth]{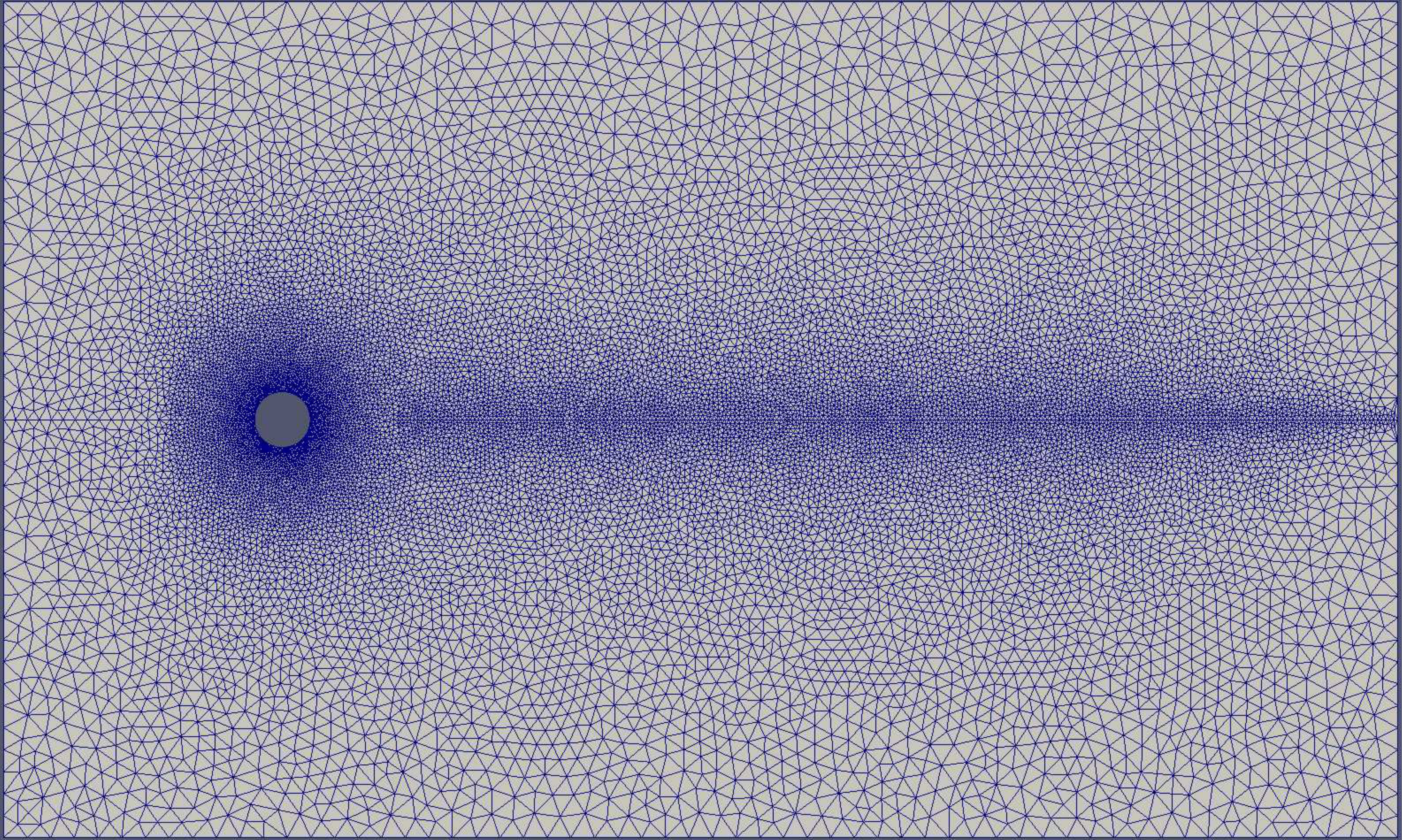}\\(a)}
	\parbox[b]{0.49\textwidth}{\centering \includegraphics[width=0.486\textwidth]{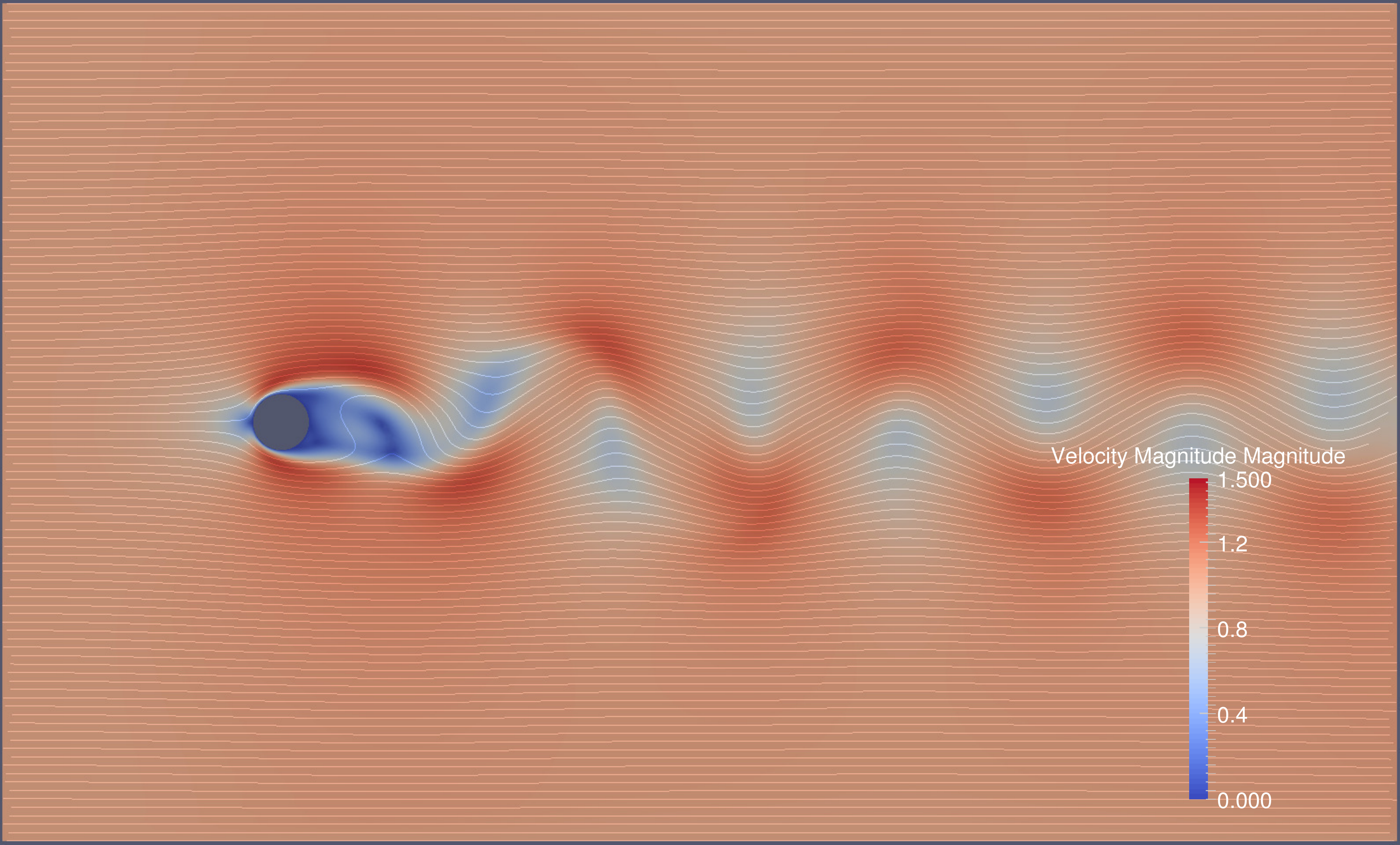}\\(b)}
	\caption{(a) FEM discretization of the domain $\Omega$ by a triangular mesh ($N = 17.048$). (b) A snapshot of the solution to \eqref{eq:NSE1} -- \eqref{eq:BCC} for a non-rotating cylinder ($\gamma(t) = 0$), the coloring is according to the velocity magnitude. The pattern is the well-known \emph{von K\'{a}rm\'{a}n vortex street}.}
	\label{fig:Mesh}
\end{figure}

The system (\ref{eq:NSE1} -- \ref{eq:BCC}) is solved by the software package \emph{OpenFOAM} \cite{JJT07} using a finite volume discretization and the \emph{PISO} scheme \cite{FP02}. We have chosen \emph{OpenFOAM} because it contains a variety of efficient solvers for various fluid flow applications. Since we will utilize finite elements for the computation of the reduced order model (evaluation of inner products etc.), we then map our solution to a finite element mesh with $N = 17.048$ degrees of freedom (Figure~\ref{fig:Mesh}(a)). This is done in the spirit of data driven modeling, where we collect data which does not necessarily have to come from a numerical method. Finally, the velocity field can be written in terms of the FEM basis:
\begin{align}
	\boldsymbol{U}(\boldsymbol{x},t) = \left( \begin{array}{c} \sum_{j=1}^N {U}^d_j(t) \phi_j(\boldsymbol{x}) \\ \sum_{j=1}^N {U}^d_{j+N}(t) \phi_j(\boldsymbol{x}) \end{array} \right), \tag{2.3} \label{eq:FEM_basis}
\end{align}
where $\left\lbrace \phi_j(\boldsymbol{x}) \right\rbrace_{j=1}^N$ are the FEM basis functions and ${U}^d(t) \in \mathbb{R}^{2N}$ are the nodal values of the two velocity components, the superscript $d$ denoting that this is a quantity defined on the grid nodes. In the following, the nodal values of all quantities will be denoted by a superscript $d$. 
All finite element related computations are performed with the \emph{FEniCS} toolbox \cite{LMW12} using linear basis functions.

For a non-rotating cylinder, i.e.~$\gamma(t) = 0$, the system possesses a periodic solution, the well-known \emph{von K\'{a}rm\'{a}n vortex street} (Figure~\ref{fig:Mesh}(b)), where vortices separate alternatingly from the upper and lower surface of the cylinder, respectively. The effect is observed frequently in nature and is one of the most studied phenomena in fluid mechanics, also in the context of flow control, where the objective is to stabilize the flow and to reduce the drag.

\subsection{Multiobjective optimal control problem}	
In many applications, the control cost is of great interest. This is immediately clear when the goal of the optimization is to save energy such that in this case, the control effort needs to be taken into account. In scalar optimization problems, this is often done by adding an additional term of the form $\beta \int_{t_0}^{t_e} \gamma^2(t) \, dt$ to the cost functional where $\beta\in\mathbb{R}_{\ge 0}$ is a weighting parameter. Here, we want to consider the two objectives flow stabilization, i.e.~the minimization of the fluctuations $\boldsymbol{u}(\boldsymbol{x},t) = \boldsymbol{U}(\boldsymbol{x},t) - \left\langle \boldsymbol{U}(\boldsymbol{x}) \right\rangle$ around the mean flow field $\left\langle \boldsymbol{U}(\boldsymbol{x}) \right\rangle = \frac{1}{T} \int_0^\top \boldsymbol{U}(\boldsymbol{x},t) \, dt$, and the minimization of the control cost separately which leads to the following multiobjective optimal control problem:
\begin{align*}
	\min_{\boldsymbol{U}, \gamma} \widehat{J}(\boldsymbol{U}, \gamma) = \min_{\boldsymbol{U}, \gamma} \left( \begin{array}{c} \int_{t_0}^{t_e} \| \boldsymbol{u}(\cdot,t)\|_{L^2}^2 \, dt \\ \|\gamma\|^2_{L^2} \end{array} \right),
\end{align*}
where $\widehat{J}: H^2(\Omega \times [t_0,t_e], \mathbb{R}^2) \times L^2([t_0,t_e],\mathbb{R}) \rightarrow \mathbb{R}^2$ and $\boldsymbol{U}(\boldsymbol{x},t)$ satisfies (\ref{eq:NSE1} -- \ref{eq:BCC}). 
In contrast to bounded domains (cf.~\cite{FGH98}), the proof of existence of a solution is an open problem for cases with \emph{no-shear} or \emph{do nothing} boundary conditions \cite{Ran00}. Nevertheless, based on numerical experiences \cite{Ran00}, we will from now on assume that there exists a unique solution $\boldsymbol{U}(\boldsymbol{x},t)$ for each $\gamma(t)$ and hence, we denote by $\boldsymbol{U}(\gamma)$ the solution $\boldsymbol{U}(\boldsymbol{x}, t)$ for a fixed $\gamma \in L^2([t_0,t_e],\mathbb{R})$ and consider the reduced cost functional $J: L^2([t_0,t_e],\mathbb{R}) \rightarrow \mathbb{R}^2$ which leads to the following multiobjective optimal control problem:
\begin{align}
	\min_{\gamma} J(\gamma) = \min_{\gamma} \left( \begin{array}{c} \int_{t_0}^{t_e} \| \boldsymbol{u}\|_{L^2}^2 \, dt \\ \|\gamma\|^2_{L^2} \end{array} \right).  \tag{MOCP} \label{eq:MOCP}
\end{align}
In general, the solution to this problem does not consist of isolated points but a set of \emph{optimal compromises} between the two objectives. In the following section, we give a short introduction to multiobjective optimal control theory and solution methods.

%% file: 03_MOCP.tex
This section is concerned with the treatment of general multiobjective optimal control problems. We will give a short introduction to the general theory before addressing the two algorithms used later on in combination with model order reduction techniques.

\subsection{Theory of multiobjective optimal control}
\label{subsec:MOCP_Theory}
\begin{figure}[h!]
	\centering
	\parbox[b]{0.3\textwidth}{\centering \includegraphics[width=0.25\textwidth]{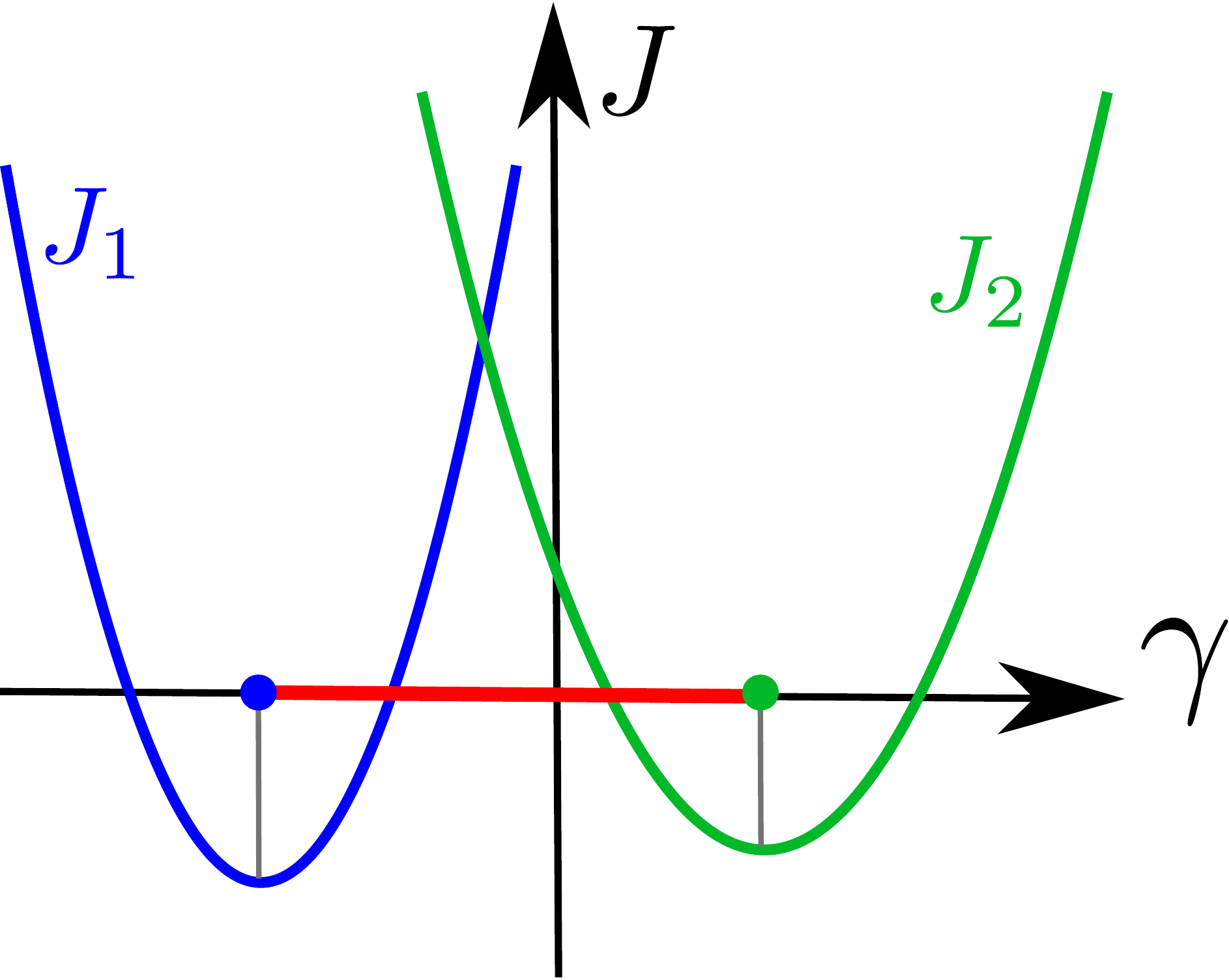}\\(a)}
	\parbox[b]{0.3\textwidth}{\centering \includegraphics[width=0.25\textwidth]{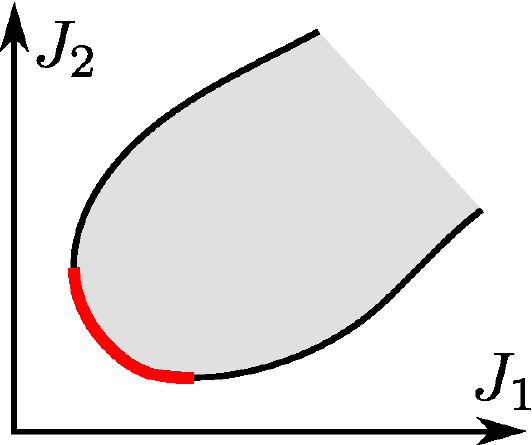}\\(b)}
	\caption{Pareto set (a) and front (b) of the multiobjective optimization problem $\min_{\gamma\in\mathbb{R}} J(\gamma)$, $J: \mathbb{R}\rightarrow\mathbb{R}^2$.}
	\label{fig:MOCP_FeasibleSet}
\end{figure}

Consider the general multiobjective optimal control problem:
\begin{align}
	\min_{\gamma} J(\gamma) = \min_{\gamma} \left( \begin{array}{c} J_1(\gamma) \\ \vdots \\ J_k(\gamma) \end{array} \right), \label{eq:generalMOCP}
\end{align}
with $J: L^2([t_0,t_e],\mathbb{R}) \rightarrow \mathbb{R}^k$ and $J_i: L^2([t_0,t_e],\mathbb{R}) \rightarrow \mathbb{R}$, $i=1,\ldots,k$. The space in which the control functions live is denoted as the \emph{decision space} and the function $J$ is a mapping to the $k$-dimensional \emph{objective space}.
The set of feasible functions $\gamma$ is the \emph{feasible set in decision space}. We denote its image as the \emph{feasible set in objective space} which consists of the \emph{feasible points} $J(\gamma)$.	
In contrast to single objective optimization problems, there exists no total order of the objective function values in $\mathbb{R}^k, k \geq 2$. Therefore, the comparison of values is defined in the following way \cite{Mie12}:
\begin{definition}~
Let $v, w \in \mathbb{R}^k$. The vector $v$ is \emph{less than} $w$ $(v <_p w)$, if $v_i < w_i$ for all $i \in \left\lbrace 1, \ldots, k \right\rbrace$. The relation $\leq_p$ is defined in an analogous way.\\
\end{definition}
A consequence of the lack of a total order is that we cannot expect to find isolated optimal points. Instead, the solution to \eqref{eq:generalMOCP} is the set of optimal compromises, the so-called \emph{Pareto set} named after Vilfredo Pareto:
\begin{definition}~
	\begin{enumerate}[(a)]
		\item A function $\gamma^*$ \emph{dominates} a function $\gamma$, if $J(\gamma^*) \leq_p J(\gamma)$ and $J(\gamma^*) \neq J(\gamma)$.
		\item A feasible function $\gamma^*$ is called \emph{(globally) Pareto optimal} if there exists no feasible function $\gamma$ dominating $\gamma^*$. The image $J(\gamma^*)$ of a (globally) Pareto optimal function $\gamma^*$ is a \emph{(globally) Pareto optimal point}.
		\item The set of nondominated feasible functions is called the \emph{Pareto set}, its image the \emph{Pareto front}.
	\end{enumerate}
\end{definition}
\noindent Consequently, for each function that is contained in the Pareto set (cf. the red line in ~Figure~\ref{fig:MOCP_FeasibleSet}(a)), one can only improve one objective by accepting a trade-off in at least one other objective. Figuratively speaking, in a two-dimensional problem, we are interested in finding the "lower left" boundary of the feasible set in objective space (cf.~Figure~\ref{fig:MOCP_FeasibleSet}(b)). More detailed introductions to multiobjective optimization can be found in \cite{Mie12, Ehr05}.

\subsection{Solution methods}
\label{subsec:MOCP_solution_methods}
Various methods exist to solve problem~\eqref{eq:generalMOCP}. In this work, we present two approaches that are fundamentally different. The first one is a reference point method \cite{RBW+09} for which the distance between a feasible point (i.e.~an objective value that lies in the feasible set in objective space) and an infeasible target in objective space is minimized. The method yields a moderate number of single objective optimization problems that are solved consecutively. The second approach is a global, derivative free subdivision algorithm \cite{DSH05} for which the Pareto set is approximated by a nested sequence of increasingly refined box coverings. However, its applicability depends critically on both the decision space dimension and the numerical effort of function evaluations. The main properties are summarized in Table~\ref{tab:MOPC_Algorithms}.
\begin{table}[htbp]
\caption{Most important properties of the two algorithms presented in this section.}
\begin{center}\footnotesize
\renewcommand{\arraystretch}{1.3}
\begin{tabular}{|c|c|c|}\hline
	~ & \textbf{Reference point method} & \textbf{Subdivision algorithm} \\
	\hline
	Optimality & local & global  \\ 
	Gradients & yes & no  \\ 
	Solution concept & Boundary of feasible set of $J(\gamma)$ & Nondominated subsets of decision space \\
	Parameter Dim. & high & moderate \\
	\hline
\end{tabular}
\end{center}
\end{table}
\label{tab:MOPC_Algorithms}

\subsubsection{Reference point method}
\label{subsubsec:MOCP_reference_point}

\begin{figure}[h!]
	\centering
	\parbox[b]{0.4\textwidth}{\centering \includegraphics[width=0.35\textwidth]{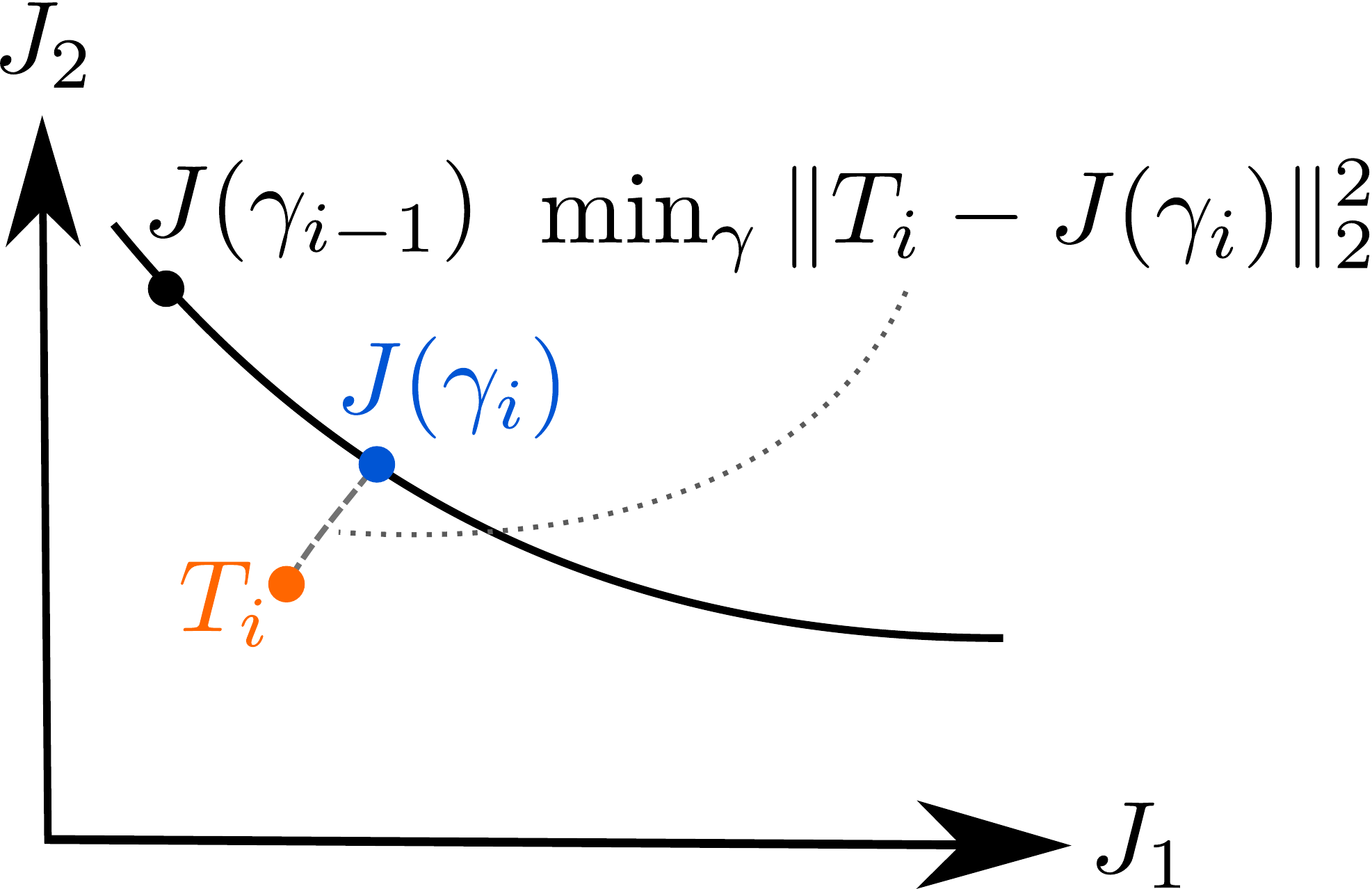}\\(a)}
	\parbox[b]{0.3\textwidth}{\centering \includegraphics[width=0.26\textwidth]{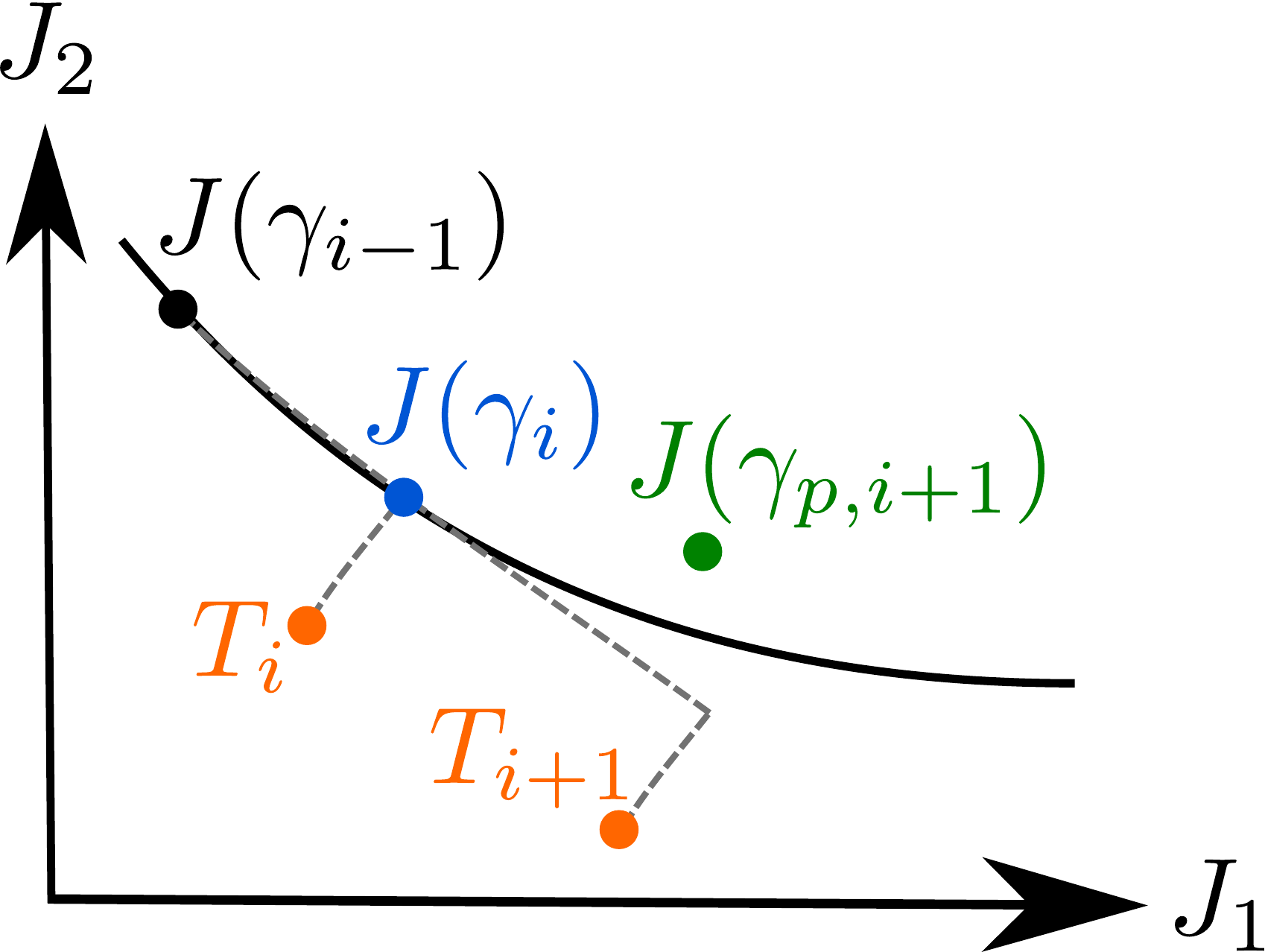}\\(b)}
	\caption{Reference point method in image space. (a) Determination of the $i$-th point on the Pareto front by solving a scalar optimization problem. (b) Computation of new target point $T_{i+1}$ and predictor step in decision space ($\gamma_{p,i+1}$, cf.~line 16 in Algorithm~\ref{Alg:ReferencePoint}).}
	\label{fig:MOCP_ReferencePoint}
\end{figure}

The reference point method presented here belongs to the category of scalarization techniques for which the solution set of \eqref{eq:generalMOCP} is approximated by a finite set of points, each computed by solving a scalar optimization problem. In the beginning, one Pareto optimal point $\gamma_0$ has to be known. This can be achieved by solving a scalar optimization problem for some weighted sum of all objectives (i.e.~$\min_{\gamma_0} s J_1(\gamma_0) + (1-s)J_2(\gamma_0),~s \in [0,1]$ for the case of two objectives), including the scalar optimization with respect to any of the objectives of \eqref{eq:MOCP}. Then, a so-called \emph{target} $T_1 \in \mathbb{R}^k$ is chosen such that it lies outside the feasible set in objective space, e.g. by shifting the solution of the first Pareto point ($T_1 = J(\gamma_0) - (h_{\parallel}, 0, \ldots, 0)^\top, h_{\parallel} > 0$). We then solve the scalar optimization problem $\min_{\gamma_{1}} \|T_{1} - J(\gamma_{1})\|^2_2$. As a result, the corresponding optimal point $J(\gamma_1)$ lies on the boundary of the feasible set and it is not possible to further improve all objectives at the same time. Thus, the point is (locally) Pareto optimal. By adjusting the target position based on targets and Pareto points already known, multiple points on the Pareto front (i.e.~$J(\gamma_2),J(\gamma_3),\ldots$) are computed recursively (cf.~Figure~\ref{fig:MOCP_ReferencePoint} for an illustration). For $J$ being continuous, the change in decision space is small when the target position changes only slightly \cite{Hil01} and hence, the current solution is a good initial guess for the next scalar problem which accelerates the convergence considerably.

Theoretically, the method is not restricted to low objective space dimensions. However, setting the targets properly to get a proper approximation of the Pareto front (i.e.~an approximation of the entire front by evenly distributed points) becomes complicated in higher dimensions. In our case, we are dealing with two objectives and the targets can be determined easily using linear extrapolation (cf.~Figure~\ref{fig:MOCP_ReferencePoint} and Algorithm~\ref{Alg:ReferencePoint}) as proposed in \cite{RBW+09}. This also allows us to compute the whole front in at most two search directions ($run=1,2$, line 2 in Algorithm~\ref{Alg:ReferencePoint}). From the initial point, we first proceed in one direction, e.g.~decreasing $J_1$. When, at some point, $J_1$ is increasing again (line 6 in Algorithm~\ref{Alg:ReferencePoint}), we have reached the \emph{extremal point} of the feasible set (cf.~Figure~\ref{fig:MOCP_FeasibleSet}(b)) and the scalar optimum of $J_1$. We then return to the initial point and proceed in the opposite direction (lines 8, 9 in Algorithm~\ref{Alg:ReferencePoint}) until the other extremal point is reached.

\begin{remark}
Using system knowledge, the algorithm can be further simplified. In the case of \eqref{eq:MOCP}, for example, we know that $\gamma(t) = 0$ is Pareto optimal since it is the scalar minimum of $J_2(\gamma)$ and therefore has the lowest possible value of $J_2$. Knowing this, we do not need to solve the initial optimization problem and moreover, we only need to move along the front in one direction until the other extremal point is reached.
\end{remark}

The scalar optimization problems can be solved using any suitable method. Here, we use a line search approach \cite{NW06} and compute the derivatives using an adjoint approach. Note that since the scalar optimization routine often is of local nature, the method overall is also local and depends on the initial guesses as well as the choice of the target points. However, under certain smoothness assumptions, the Pareto set is a $(k-1)$-dimensional manifold \cite{Hil01} such that once the first point is computed correctly, the method is promising to find the globally optimal Pareto set for many problems.

\begin{remark}
	For the reference point method, it turned out to be numerically beneficial that $(T_i - J_i) / (T_j - J_j)\ = {\mathcal O}(1) \ \forall \ i,j \in \lbrace 1, ..., k \rbrace$. Otherwise, the computation of new target points may become sensitive to the step length parameters.
\label{remark:ReferencePoint}
\end{remark}

\begin{algorithm}[H]
\caption{(Reference point method for $J \in \mathbb{R}^2$)}
\label{Alg:ReferencePoint}
\begin{algorithmic}[1]
\Require Initial solution $\gamma_0$, step length parameters $h_{\parallel}, h_{\perp}, h_{p}>0$, index $i = 0$
\State Compute the first target point $T_1 = J(\gamma_0) - (h_{\parallel}, 0)^\top$.
\For{$run=1, 2$}
\Loop
\State $i = i + 1$
\State Solve scalar optimization problem $\min_{\gamma_{i}} \|T_{i} - J(\gamma_{i})\|^2_2$
\If{extremal point of Pareto front is passed ($J_{run}(\gamma_i) > J_{run}(\gamma_{i-1})$)}
\If{$run=1$ (first direction is completed)}
\State $\gamma_{p,i+1} = \gamma_0$ (Go back to the initial solution)
\State $T_{i+1} = J(\gamma_0) - h_{\parallel} \frac{J(\gamma_1) - J(\gamma_{0})}{\|J(\gamma_1) - J(\gamma_{0})\|_2} + h_{\perp} \frac{T_1 - J(\gamma_{1})}{\|T_1 - J(\gamma_{1})\|_2}$ (Opposite direction)
\State \textbf{break}
\Else
\State \textbf{STOP}
\EndIf
\Else
\State $T_{i+1} = J(\gamma_i) + h_{\parallel} \frac{J(\gamma_i) - J(\gamma_{i-1})}{\|J(\gamma_i) - J(\gamma_{i-1})\|_2} + h_{\perp} \frac{T_i - J(\gamma_{i})}{\|T_i - J(\gamma_{i})\|_2}$
\State $\gamma_{p,i+1} = \gamma_i + h_{p} \left(\gamma_i - \gamma_{i-1}\right)$ (Predictor step)
\EndIf
\EndLoop
\EndFor
\end{algorithmic}
\end{algorithm}

\subsubsection{Global subdivision algorithm}
\label{subsubsec:MOCP_subdivision}
When the image space dimension increases or, alternatively, the Pareto front is disconnected (see e.g.~\cite{DEF+16}), continuation methods like the reference point method may fail. Moreover, the resulting scalar optimization problems are often solved by algorithms of local nature which may also not be sufficient if global optima are desired. In addition to that, derivatives are hard to compute or not available at all in many applications. The subdivision algorithm presented here overcomes these problems. It is described in more detail in \cite{DSH05} including a proof of convergence. The version using gradient information is based on concepts for the computation of attractors of dynamical systems \cite{DH97}, it is however not considered here. Instead, we focus on the derivative free approach. It is very robust and since we directly utilize the concept of dominance when comparing solutions, we theoretically do not need to make any assumptions about the problem except that the control is real-valued, i.e.  $\gamma \in \mathbb{R}^m$. In the numerical realization however, we approximate sets by sample points and hence, we further require that the objective $J$ is a continuous function of $\gamma$.

\begin{figure}[h!]
	\centering
	\parbox[b]{0.24\textwidth}{\centering \includegraphics[width=0.24\textwidth]{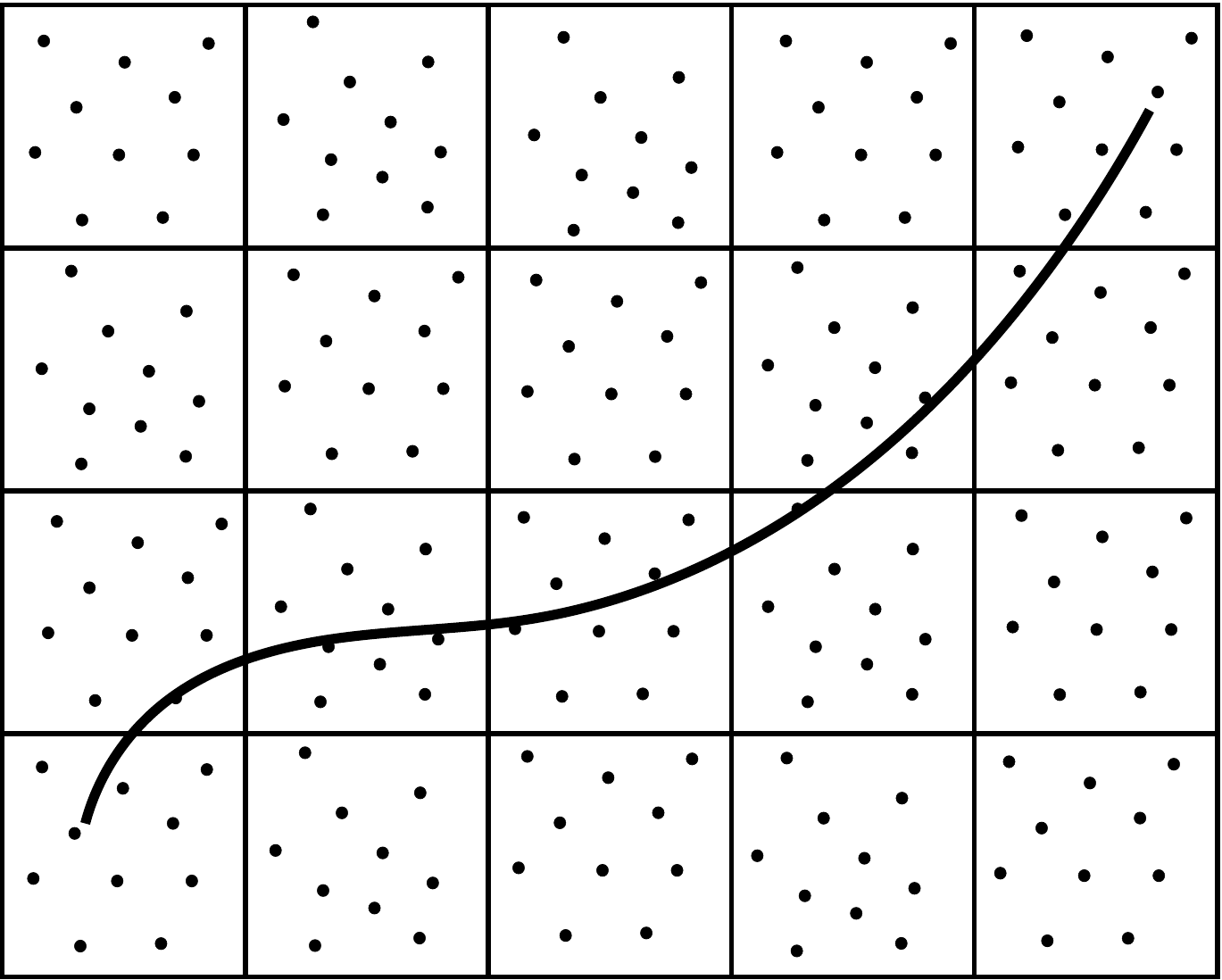}\\(a)}
	\parbox[b]{0.24\textwidth}{\centering \includegraphics[width=0.24\textwidth]{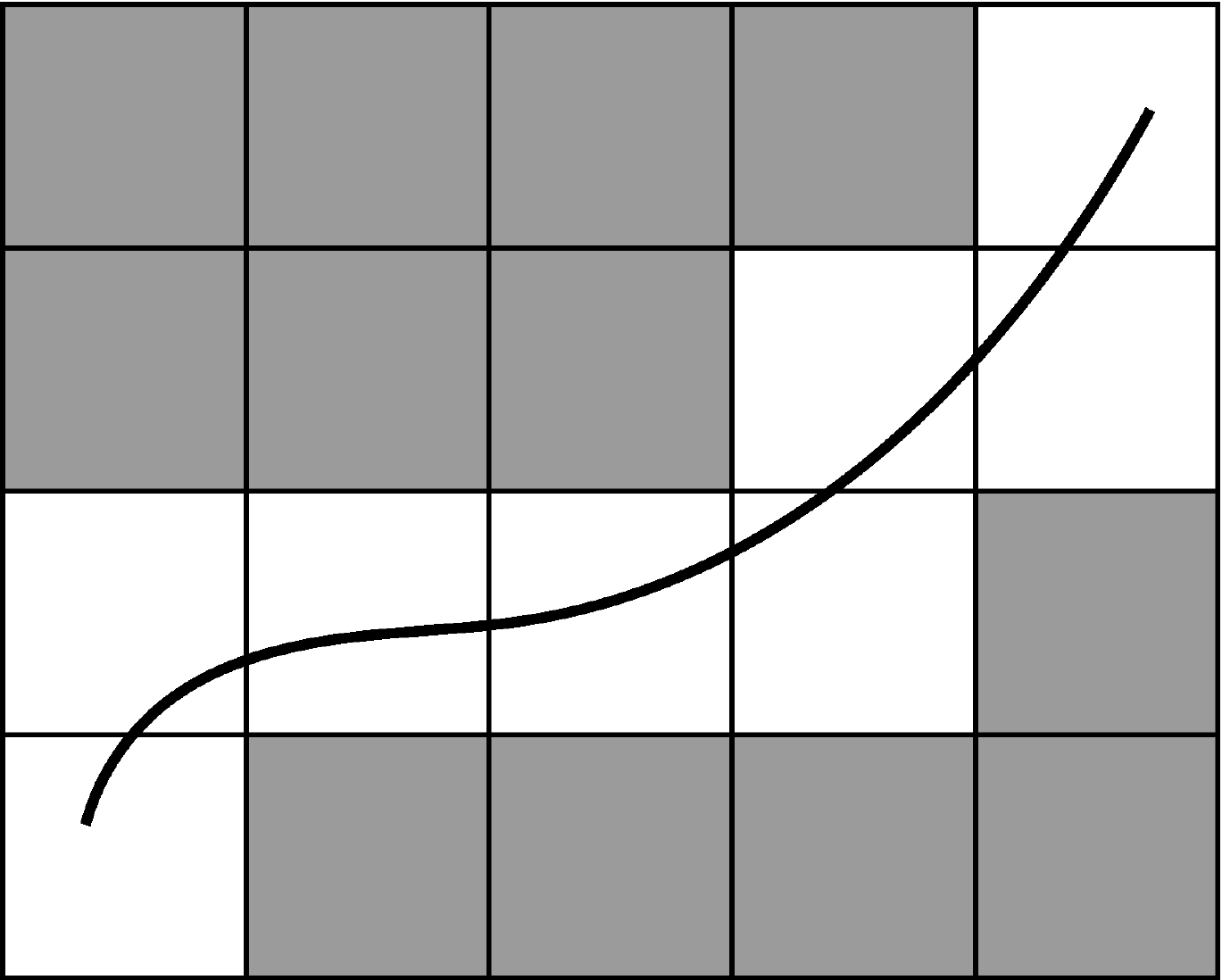}\\(b)}
	\parbox[b]{0.24\textwidth}{\centering \includegraphics[width=0.24\textwidth]{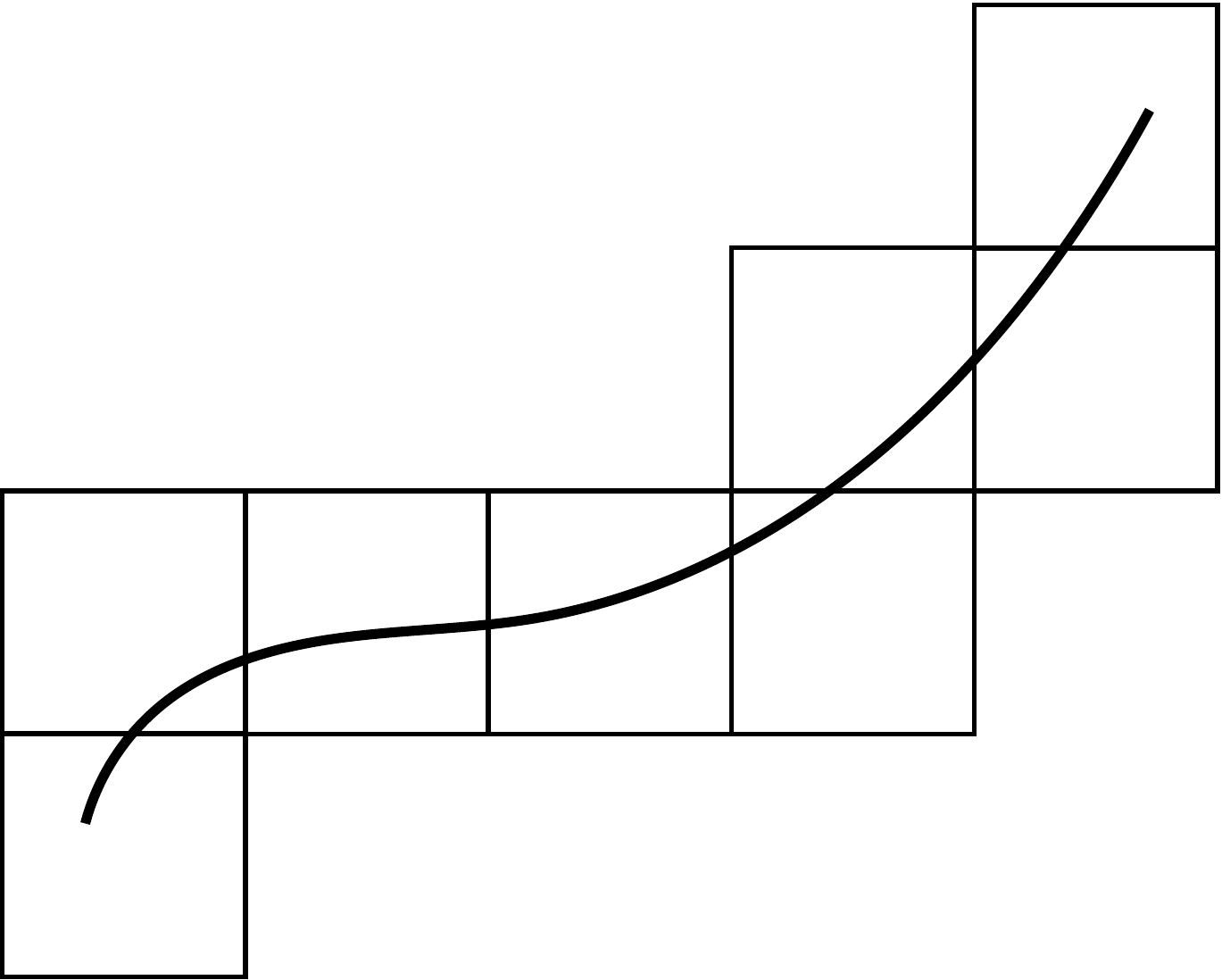}\\(c)}
	\parbox[b]{0.24\textwidth}{\centering \includegraphics[width=0.24\textwidth]{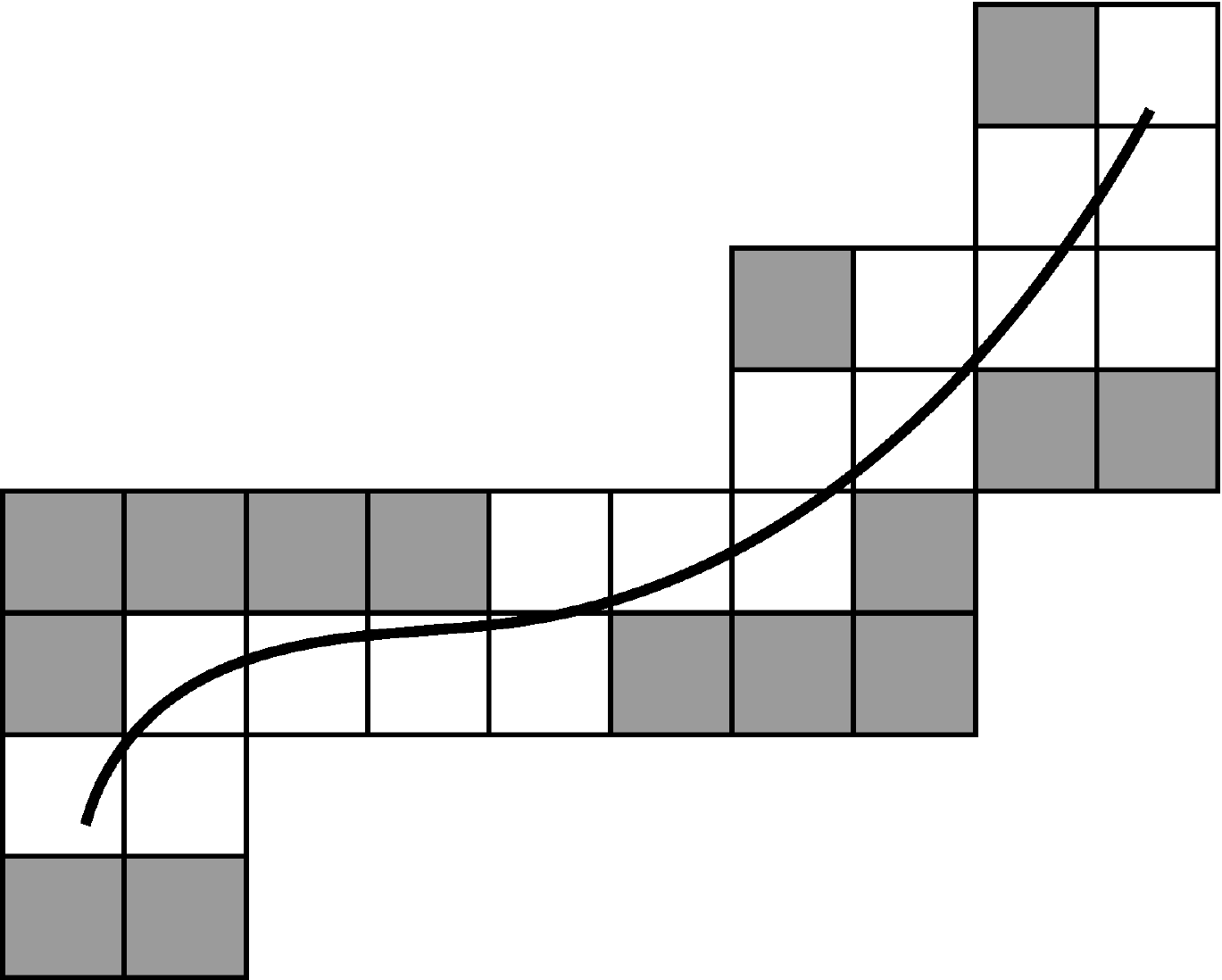}\\(d)}
	\caption{Global subdivision algorithm, example with $\gamma \in \mathbb{R}^2$. (a) Sampling. (b) Nondominance test. (c) Elimination of dominated boxes. (d) Subdivision, sampling and nondominance test.}
	\label{fig:MOCP_GAIO}
\end{figure}

The subdivision algorithm (cf.~Algorithm~\ref{Alg:GAIO} and Figure~\ref{fig:MOCP_GAIO}) is a set oriented method where the Pareto set is approximated by a nested sequence of increasingly refined box coverings (cf.~Figure~\ref{fig:MOCP_GAIO}). During the algorithm, a sequence of box collections $\{\mathcal{B}_i\}_{i=0,\ldots}$ covering the Pareto set is constructed, starting with a sufficiently large initial set $\mathcal{B}_0\subset\mathbb{R}^m$. (This results in box constraints for the control $\gamma$.) In each iteration, the algorithm performs the steps \emph{subdivision} and \emph{selection} until a given precision criterion is satisfied. This way, a subset of the previous box collection remains that is a closer covering of the desired set. In the selection step, the boxes are compared pairwise and all boxes that are dominated by another box are eliminated from the collection. Numerically, this is realized by a representative set of sampling points in each box (cf.~Figure~\ref{fig:MOCP_GAIO}(a)), at which the objectives $J_i$ are evaluated. In this case, we say that a box is dominated if all sampling points are dominated by at least one point from another box. The remaining boxes are then subdivided into two boxes of half the size and we proceed with the next selection step. The subdivision step is performed cyclically in the decision space dimensions $1$ to $m$. Hence, in order to divide each dimension $q$ times, we require $N_{sub} = mq$ subdivision steps in total.
\begin{algorithm}
	\caption{(Global, derivative free subdivision algorithm)}
	\label{Alg:GAIO}
	\begin{algorithmic}[1]
	\Require Box constraints $\gamma^{min}, \gamma^{max} \in \mathbb{R}^m$, number of subdivision steps $N_{sub}$
	\State Create initial box collection $\mathcal{B}_0$ defined by the constraints $\gamma^{min}, \gamma^{max}$, i.e.~$\mathcal{B}_0 = B = [\gamma^{min}_1, \gamma^{max}_1] \times \ldots \times [\gamma^{min}_m, \gamma^{max}_m]$
	\For{$i=1,\ldots,N_{sub}$}
	\State Subdivision step: $\bigcup_{B \in \widehat{\mathcal{B}}_{i}} B = \bigcup_{B \in \mathcal{B}_{i-1}} B$ such that $\max_{B\in\widehat{\mathcal{B}}_i}(\text{diam}(B)) < \Theta \max_{B\in\mathcal{B}_{i-1}}(\text{diam}(B)), \ \Theta \in (0,1)$
	\State Insert $S$ sampling points $\gamma_1,\ldots,\gamma_S$ each box and evaluate the cost functional $J(\gamma_s), s = 1,\ldots,S$
	\State Selection step: Eliminate all boxes that contain only dominated points: $\mathcal{B}_i = \left\lbrace B \in \widehat{\mathcal{B}}_i \big| \text{There exists no } \gamma^* \in \bigcup_{\widehat{B} \in \widehat{\mathcal{B}}_{i} \setminus B} \widehat{B} \text{ such that } J(\gamma^*) \leq_p J(\gamma)~\forall \gamma \in B \right\rbrace$
	\EndFor
	\end{algorithmic}
\end{algorithm}

Since the applicability of the subdivision algorithm depends critically on both the decision space dimension $m$ and the numerical effort of function evaluations, 
it is in practice limited to low decision space dimensions. Therefore, we introduce a sinusoidal control, i.e.~$\gamma(t) = A\sin(2 \pi \omega t + \tau)$. The choice is motivated by the fact that the uncontrolled dynamics of the problem at hand are also periodic and we include a phase shift to allow for controls with non-zero initial conditions. This way, we transform \eqref{eq:MOCP} into an MOP with $\widetilde{J}: \mathbb{R}^3 \rightarrow \mathbb{R}^2$:
\begin{align}
	\min_{A, \omega, \tau \in \mathbb{R}} \widetilde{J}(A, \omega, \tau) = \min_{A, \omega, \tau \in \mathbb{R}} \left( \begin{array}{c} \int_{t_0}^{t_e} \| \boldsymbol{u}\|_{L^2}^2\ dt \\ \int_{t_0}^{t_e} \left( A\sin(2 \pi \omega t + \tau) \right)^2\ dt \end{array} \right). \tag{MOP-3D} \label{eq:MOP-2D}
\end{align}

%% file: 04_ROM.tex
The problem (\ref{eq:MOCP}) could now be solved using either of the methods presented in Section~\ref{subsec:MOCP_solution_methods}. However, both require a large number of evaluations of the cost functional and consequently, evaluations of the system (\ref{eq:NSE1} -- \ref{eq:BCC}). A finite element or finite volume discretization yields a large number of degrees of freedom such that solving \eqref{eq:MOCP} quickly becomes computationally infeasible. Hence, we need to reduce the cost for solving the dynamical system significantly. This is achieved by means of reduced order modeling where we replace the state equation by a reduced order model of coupled, nonlinear ordinary differential equations that can be solved much faster.

\subsection{Reduced Order Model}
\label{subsec:ROM}
The standard procedure for deriving a reduced order model is by introducing a Galerkin ansatz \cite{HLB98}:
\begin{align}
	\boldsymbol{U}(\boldsymbol{x},t) = \sum_{j=1}^{l} \alpha_j(t) \boldsymbol{\psi}_j(\boldsymbol{x}), \hspace{0.5cm} (\boldsymbol{x},t) \in \Omega \times [t_0,t_e], \label{eq:Galerkin}
\end{align}
where $\alpha(t)$ are time-dependent coefficients and $\boldsymbol{\psi}_j(\boldsymbol{x})$ are basis functions. In contrast to FEM, these basis functions may in general have full support, i.e.~they are global. The objective is now to find a basis as small as possible ($l \ll N$) while allowing for a good approximation of the flow field. 

Using \eqref{eq:Galerkin}, we can derive a reduced model of dimension $l$. By selecting a basis that is divergence-free, the mass conservation equation \eqref{eq:NSE2} is automatically satisfied. (When computing this basis using data from a divergence-free flow field, the resulting basis elements are also divergence-free \cite{BHL93}.) Inserting \eqref{eq:Galerkin} into the weak formulation of the Navier-Stokes equations \eqref{eq:NSE_weak} yields a system of $l$ ordinary differential equations of the form	
\begin{align*}
	\dot \alpha(t) = \widetilde{\mathcal{B}} \alpha(t) + \widetilde{\mathcal{C}} \left(\alpha(t)\right),
\end{align*}
where the coefficient matrices are computed by evaluating the $L^2$ inner products:
\begin{align}
	\left( \widetilde{\mathcal{B}} \right)_{ij} &= \frac{1}{Re} \left( \nabla \boldsymbol{\psi}_i, \nabla \boldsymbol{\psi}_j \right)_{L^2}, \notag \\
 	\left( \widetilde{\mathcal{C}} \left(\alpha(t)\right) \right)_j &= \alpha(t)^\top \mathcal{Q}_j \alpha(t), \notag \\
 	\text{with   } \left( \mathcal{Q}_j \right)_{ik} &= \left( \left( \boldsymbol{\psi}_i \cdot \nabla \right) \boldsymbol{\psi}_k, \boldsymbol{\psi}_j \right)_{L^2}. \notag
\end{align}
\begin{remark}
In the above as well as all following formulations, we have neglected the influence of the pressure term. In \cite{NPM05}, it was shown that including the pressure term is favorable for the accuracy of the reduced model for open systems. Since we consequently perform a model stabilization in order to increase the accordance with the flow field data, this results in a further modification of the model coefficients. Consequently, we have (following \cite{CEMF09}) decided to neglect the influence of the pressure term in this work.
\end{remark}

This model, however, is not controllable in the sense that we can no longer influence the solution by choosing $\gamma(t)$. Moreover, the basis functions need to be computed from data with homogeneous boundary conditions ($\boldsymbol{U}(\boldsymbol{x},t) = 0$ for $\boldsymbol{x} \in \Gamma$) such that they also possess homogeneous boundary conditions \cite{Sir87}. Otherwise, it cannot be guaranteed that the Galerkin ansatz \eqref{eq:Galerkin} satisfies the PDE boundary conditions. The idea is therefore to introduce the Galerkin expansion for a modified flow field $\boldsymbol{u}(\boldsymbol{x}, t)$ which satisfies homogeneous boundary conditions for all $t \in [t_0,t_e]$. In order to realize this, we introduce the following flow field decomposition \cite{Rav00, BCB05}:
\begin{align}
	\boldsymbol{U}(\boldsymbol{x}, t) &= \left\langle \boldsymbol{U}(\boldsymbol{x}, t) \right\rangle + \gamma(t) \boldsymbol{U}_{c}(\boldsymbol{x}) + \boldsymbol{u}(\boldsymbol{x}, t), \label{eq:extGalerkin}
\end{align}
where $\boldsymbol{U}_c(\boldsymbol{x})$ is a so-called \emph{control function} which is a solution to \eqref{eq:NSE1}, \eqref{eq:NSE2}, \eqref{eq:NSE_IC} with a constant cylinder rotation $\gamma_c$ and zero boundary conditions elsewhere. The choice of the rotating velocity for the computation of $\boldsymbol{U}_c(\boldsymbol{x})$ is arbitrary but has an influence on the resulting reduced model such that the model performance might be influenced by this choice. Since the main intention of this article is to develop multiobjective optimal control algorithms for PDE-constrained problems, we do not address this problem further and choose $\gamma(t) = \gamma_c = 2$ such that the cylinder surface rotates with a velocity of $1$. The decomposition is now computed in two steps:
\begin{align*}
	\widetilde{\boldsymbol{U}}(\boldsymbol{x}, t) &= \boldsymbol{U}(\boldsymbol{x}, t)\big|_{\gamma(t)=\gamma_{ref}(t)} - \frac{\gamma_{ref}(t)}{\gamma_c} \boldsymbol{U}_c(\boldsymbol{x}), \\
	\boldsymbol{u}(\boldsymbol{x},t) &= \widetilde{\boldsymbol{U}}(\boldsymbol{x}, t) - \left\langle \widetilde{\boldsymbol{U}}(\boldsymbol{x}, t) \right\rangle,
\end{align*}
where $\gamma_{ref}$ is the \emph{reference control} for which the data was collected.


In this way, $\boldsymbol{u}(\boldsymbol{x}, t) = \sum_{j=1}^{l} \alpha_j(t) \boldsymbol{\psi}_j(\boldsymbol{x})$ satisfies homogenous boundary conditions for all $t \in [t_0,t_e]$. (Note that this does not hold exactly on $\Gamma_N$, however, the resulting error is small and hence neglected \cite{BCB05}.) Inserting \eqref{eq:extGalerkin} into \eqref{eq:NSE_weak} yields an extended reduced order model:
\begin{align}
 	\dot \alpha(t) &= \mathcal{A} + \mathcal{B} \alpha(t) + \mathcal{C} \left(\alpha(t)\right) + \mathcal{D} \dot \gamma(t) + \left(\mathcal{E} + \mathcal{F} \alpha(t) \right) \gamma(t) + \mathcal{G} \gamma^2(t), \label{eq:eROM_state} \\
	\alpha_j(t_0) &= \left( u(\cdot, t_0), \boldsymbol{\psi}_j\right)_{L^2} =: \alpha_{0,j}, \label{eq:eROM_initial}
\end{align}
where $\alpha \in H^1([t_0,t_e],\mathbb{R}^l)$. The expressions for the coefficient matrices $\mathcal{A}$ to $\mathcal{G}$ are given in Appendix~\ref{appendix:ROM}, their size depends on the dimension $l$ of the ROM: $\mathcal{A},\ \mathcal{D},\ \mathcal{E},\ \mathcal{G} \in \mathbb{R}^l;\ \mathcal{B},\ \mathcal{F} \in \mathbb{R}^{l,l};\ (\mathcal{C} (\alpha(t)))_j = \alpha(t)^\top \mathcal{Q}_j \alpha(t)$ with $\mathcal{Q} \in \mathbb{R}^{l,l,l}$.

\begin{remark}
	Since the time derivative of the control occurs in \eqref{eq:eROM_state}, we require higher regularity, i.e.~$\gamma \in H^1([t_0,t_e],\mathbb{R})$.
	\label{remark:derivativeSpace}
\end{remark}

Due to the truncation of the POD basis, the higher modes covering the small scale dynamics, i.e.~the dynamics responsible for dissipation, are neglected. This can lead to incorrect long term behavior and hence, the model needs to be modified to account for this. Several researchers have addressed this problem and proposed different solutions, see e.g.~\cite{Rem00, SK04, NPM05, CEMF09, BBI09}. Here, we achieve stabilization of the model \eqref{eq:eROM_state} by solving an augmented optimization problem \cite{GBZI04}. In this approach the entries of the matrix $\mathcal{B}$ are manipulated in order to minimize the difference between the ROM solution and the projection of the data onto the reduced basis ($\alpha_j^p(t) = (U(\cdot,t),\boldsymbol{\psi}_j)_{L^2}$). The resulting initial value problem is solved using a Runge-Kutta method of fourth order. 

Finally, we replace the objectives in \eqref{eq:MOCP} by their equivalents in the reduced order model. 
Since the reference point method is more stable for objectives of the same order of magnitude (cf.~Remark~\ref{remark:ReferencePoint}), we additionally multiply the second component by $l$ (i.e.~the number of basis functions)
\begin{align}
	\min_{\gamma} \left( \begin{array}{c} \int_{t_0}^{t_e} \int_{\Omega} \| \sum_{i=1}^l \alpha_i(t) \boldsymbol{\psi}_i(x) \|^2_2 \, dx\, dt \\ l\int_{t_0}^{t_e} \gamma^2(t) \, dt \end{array} \right).  \label{eq:R-MOCP1} 
\end{align}

\subsection{Proper orthogonal decomposition (POD)}
\label{subsec:POD}
In order to achieve high numerical efficiency, we want to compute a basis as small as possible. To this end, we utilize the well-known Proper Orthogonal Decomposition (POD). Using POD, we can compute the 
orthonormal basis $\left\lbrace\boldsymbol{\psi}_i(\boldsymbol{x})\right\rbrace_{i=1}^l$ which, for every $l$, optimally represents a set of flow field realizations $\left\lbrace \boldsymbol{u}(\boldsymbol{x}, t_1), \dots, \boldsymbol{u}(\boldsymbol{x}, t_{m}) \right\rbrace$ \cite{Sir87}. Basically, the computation is realized via a singular value decomposition of the so-called \emph{snapshot matrix} $S \in \mathbb{R}^{2N,m}, S = \left[{u}^d(t_1), \dots, {u}^d(t_m)\right]$, where we collect $2N$ measurements (i.e.~the velocity components at $N$ nodes) at $m$ different time instants. If $m \ll N$, as it is often the case for FEM discretizations, it is more efficient to solve the $m$-dimensional eigenvalue problem \cite{KV02, Fah00}:
\begin{align}
	S^\top M S v_i = \sigma_i v_i, \hspace{1cm} i=1,\ldots,m, \label{eq:POD_Eig}
\end{align}
where $M \in \mathbb{R}^{2N,2N}$ is the finite element mass matrix. Using the eigenvalues and eigenvectors from \eqref{eq:POD_Eig}, the POD modes can be computed by making use of the FEM basis functions $\left\lbrace \phi_j(\boldsymbol{x}) \right\rbrace_{j=1}^N$:
\begin{align*}
	{\psi}^d_i &= \frac{1}{\sqrt{\sigma_i}} S v_i, \\
	\boldsymbol{\psi}_{i}(\boldsymbol{x}) &= \left( \begin{array}{c}
	\sum_{j=1}^N {\psi}^d_{i,j} \phi_j(\boldsymbol{x}) \\
	\sum_{j=1}^N {\psi}^d_{i,j+N} \phi_j(\boldsymbol{x})
	\end{array} \right).
\end{align*}

The eigenvalues $\sigma$ are a measure for the information contained in the respective modes and hence, the error between the Galerkin ansatz and the snapshot data is determined by the ratio of the truncated eigenvalues \cite{Sir87}, i.e.~$\epsilon(l) = \sum_{j=1}^l \sigma_i / \sum_{j=1}^{m} \sigma_i$. By choosing a value for $\epsilon$, typically $0.99$ or $0.999$, the basis size can be determined. For many applications, the eigenvalues decay fast such that a truncation to a small basis is possible. 

\begin{figure}[h!]
	\centering
	\includegraphics[width=0.35\textwidth]{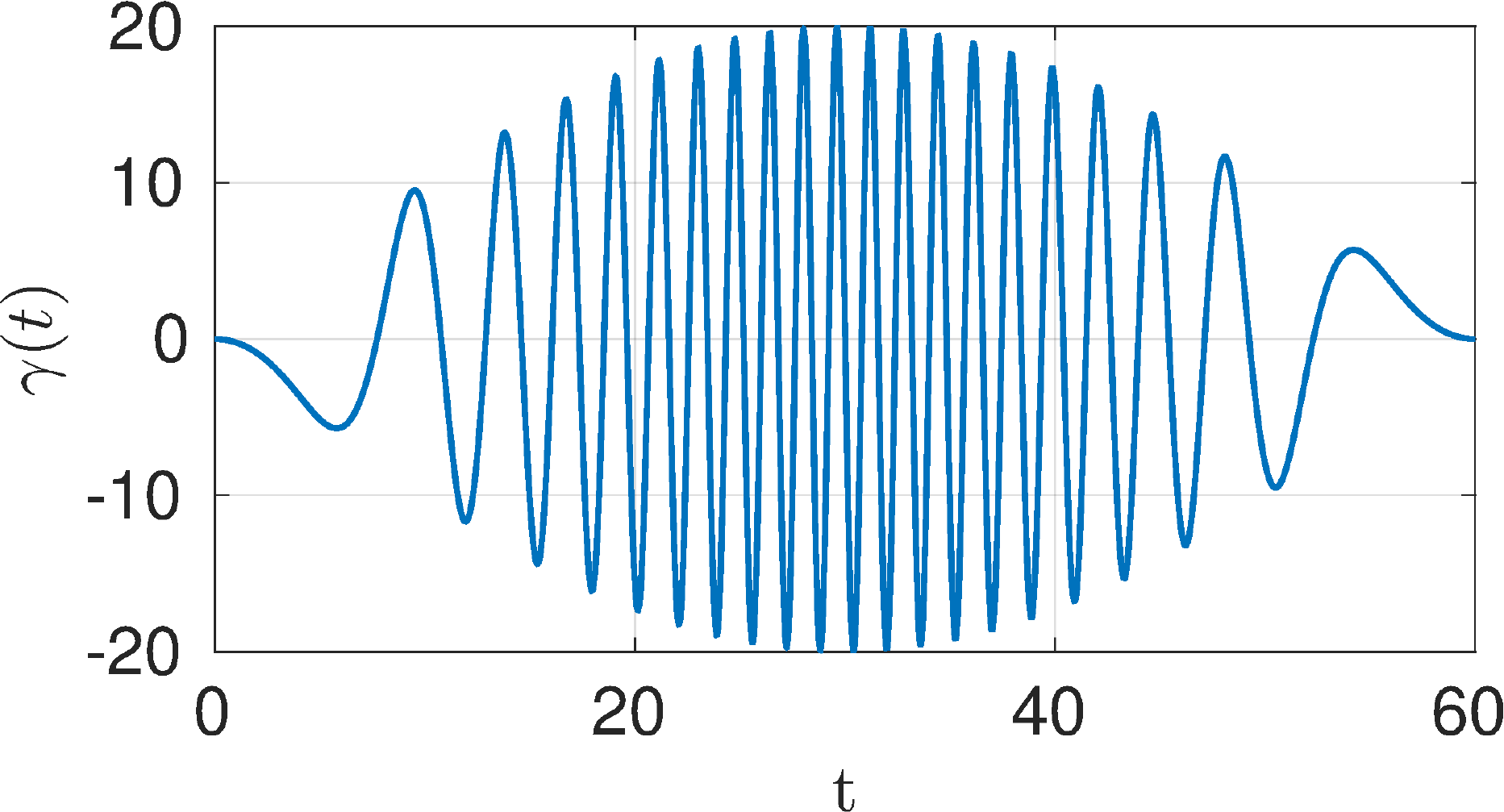}
	\caption{Chirping function ($\gamma_{chirp} = -4 \sin(2\pi t / 120) \cos( 2 \pi t / 3 - 18 \sin(2 \pi t / 60))$) used for the computation of the snapshot matrix \cite{BCB05}.}
	\label{fig:Chirp}
\end{figure}

\begin{figure}[h!]
	\centering
	\parbox[b]{0.4\textwidth}{\centering \includegraphics[width=0.35\textwidth]{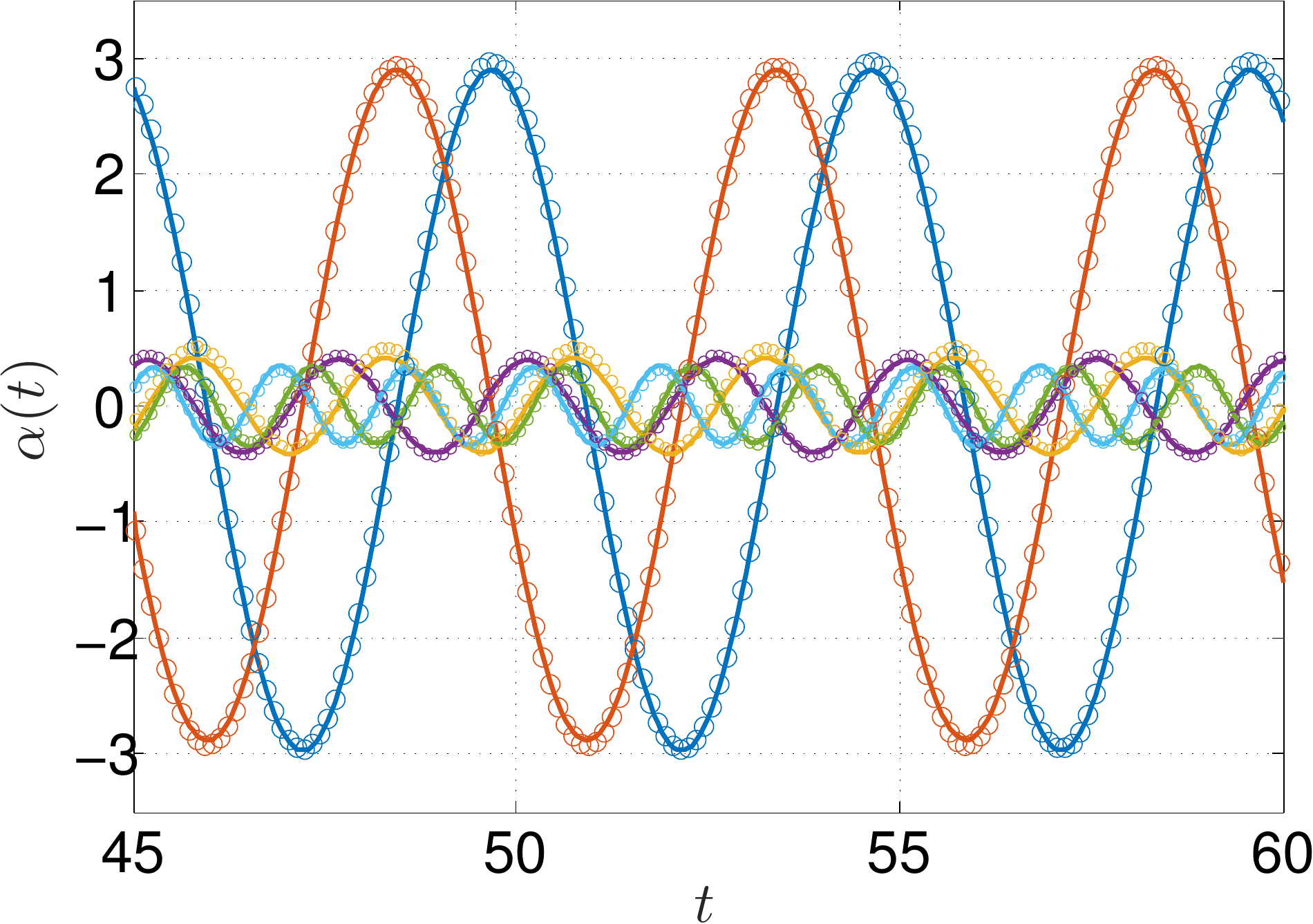}\\(a)}
	\parbox[b]{0.4\textwidth}{\centering \includegraphics[width=0.35\textwidth]{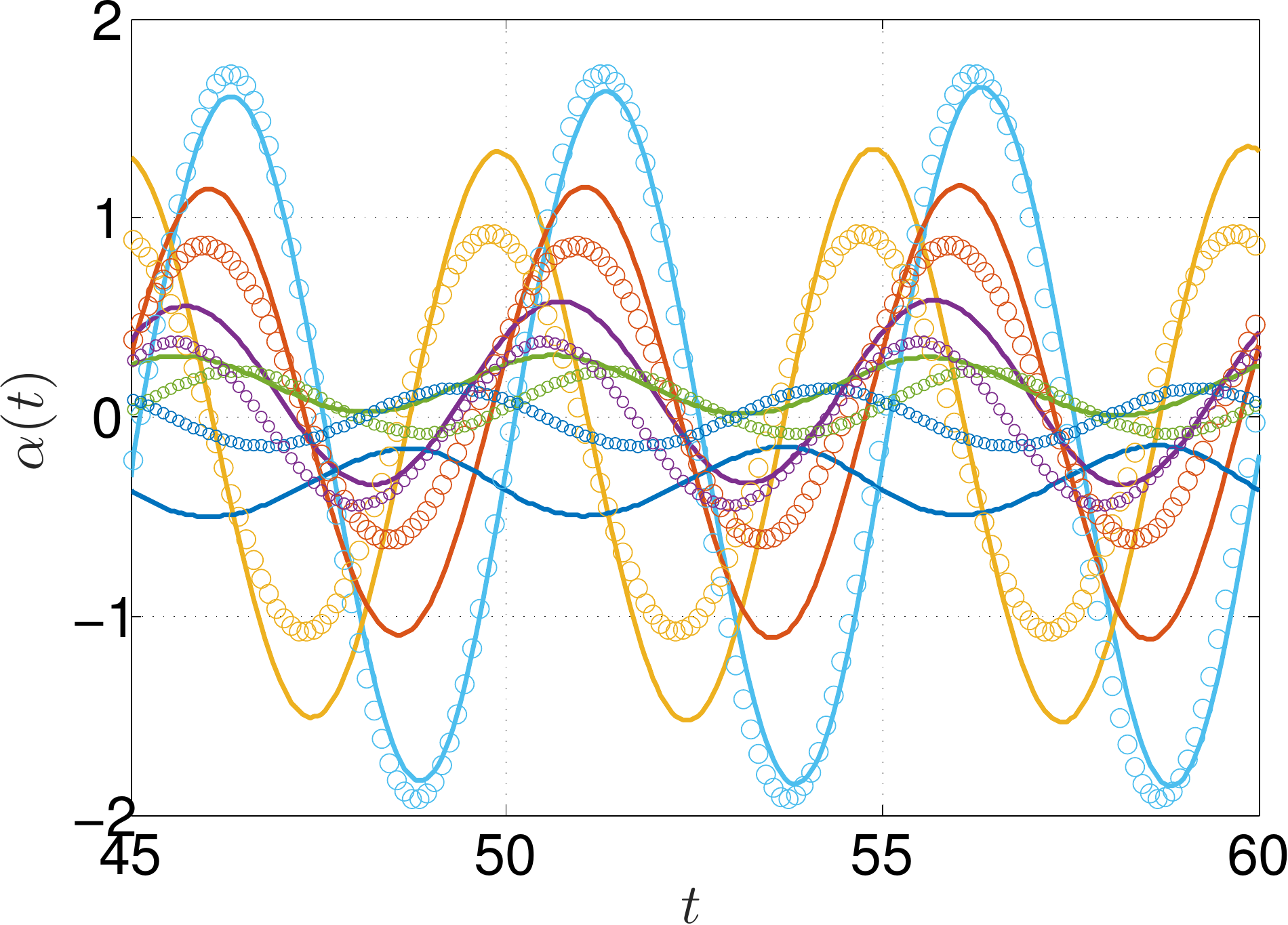}\\(b)}
	\caption{First six modes of \eqref{eq:eROM_state}, \eqref{eq:eROM_initial} (---) and projection of FEM solution ($\circ$) for two reduced order models with different reference controls. (a) $\gamma_{ref} = 0$. (b) $\gamma_{ref} = \gamma_{chirp}$.}
	\label{fig:alpha_ROM_vs_FEM}
\end{figure}

The error $\epsilon$ is only known for the reference control at which the data was collected. If one allows for multiple evaluations of the finite element solution, the reduced model can be updated when necessary. This is realized using, e.g., trust-region methods \cite{Fah00, BC08}. An alternative approach is to use a reference control that yields sufficiently rich dynamics such that the model can be expected to be valid for a larger variety of controls \cite{BCB05}. Since it is computationally cheaper, we follow the second approach and take $1201$ snapshots in the interval $[0,60]$ with $\Delta t = 0.05$ for a \emph{chirping} reference control (cf.~Figure~\ref{fig:Chirp}). This leads to a basis of size $l=38$ for $\epsilon(l) \geq 99\%$. Figure~\ref{fig:alpha_ROM_vs_FEM} depicts the comparison between the uncontrolled solution (i.e.~$\gamma(t) = 0$) of the reduced model and the projection of the finite volume solution for reduced models computed with a reference control $\gamma_{ref} = 0$ and $\gamma_{ref} = \gamma_{chirp}$, respectively.

\begin{figure}[h!]
	\centering
	\includegraphics[width=0.4\textwidth]{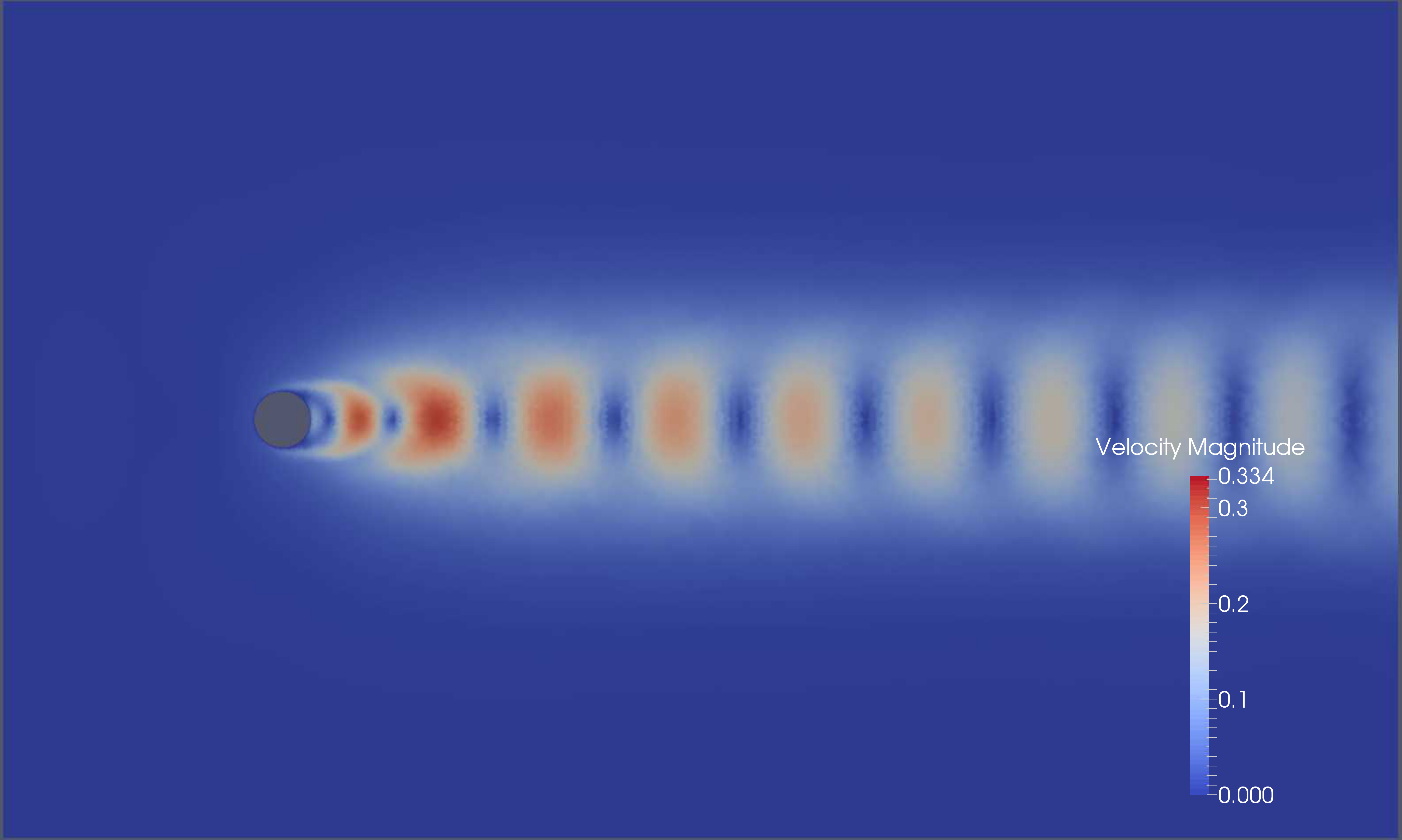}
	\includegraphics[width=0.4\textwidth]{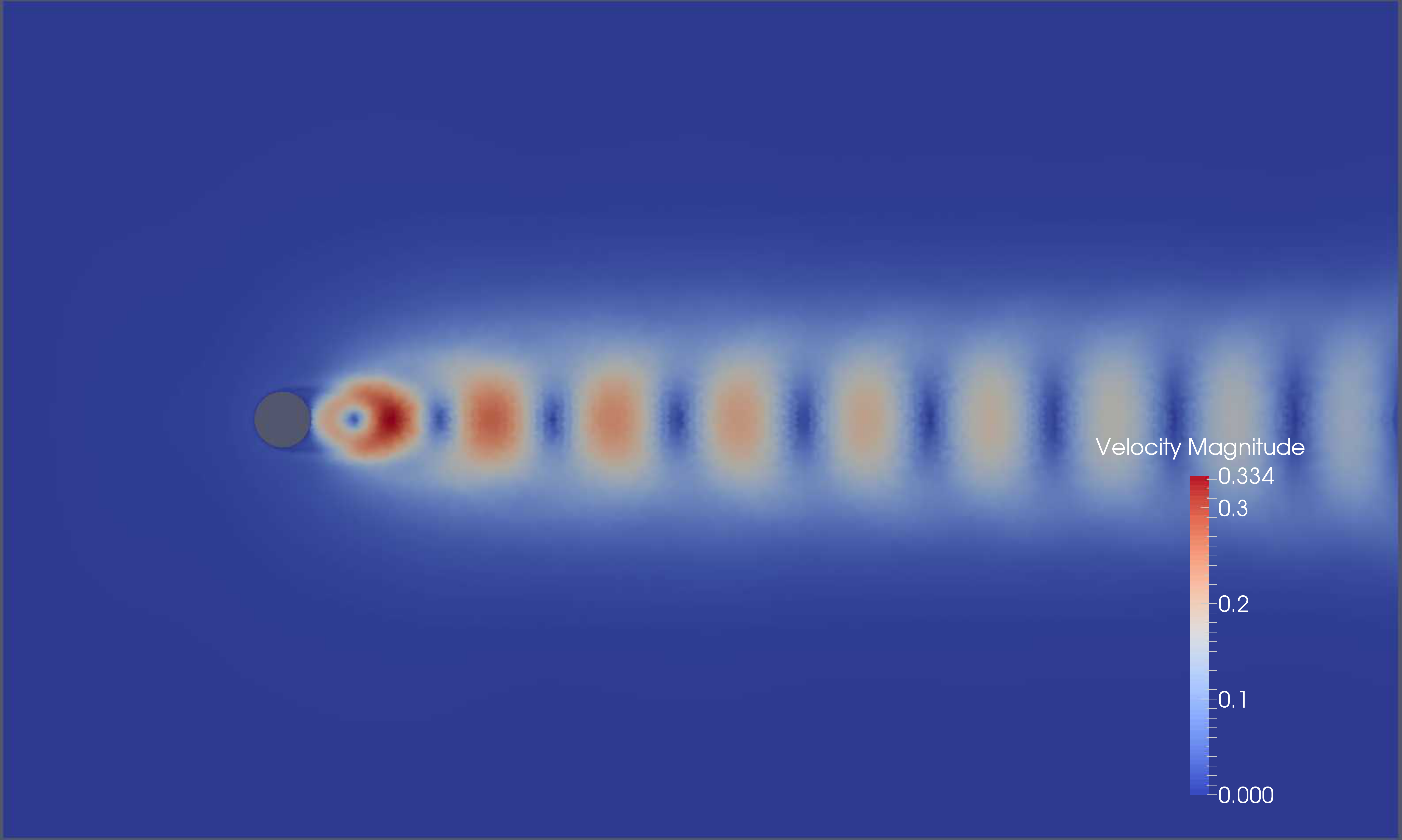}\\ \vspace{0.1cm}
	\includegraphics[width=0.4\textwidth]{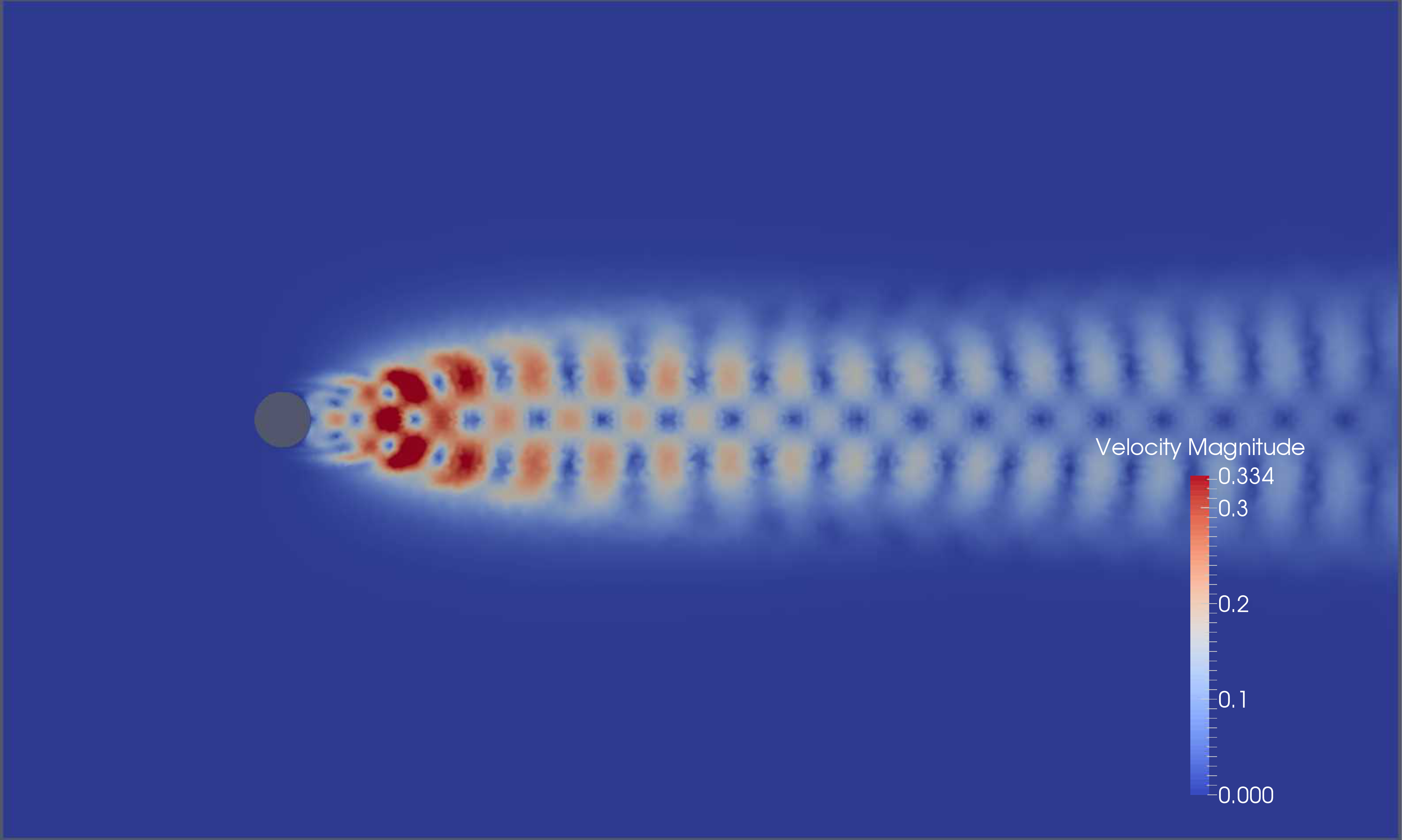}
	\includegraphics[width=0.4\textwidth]{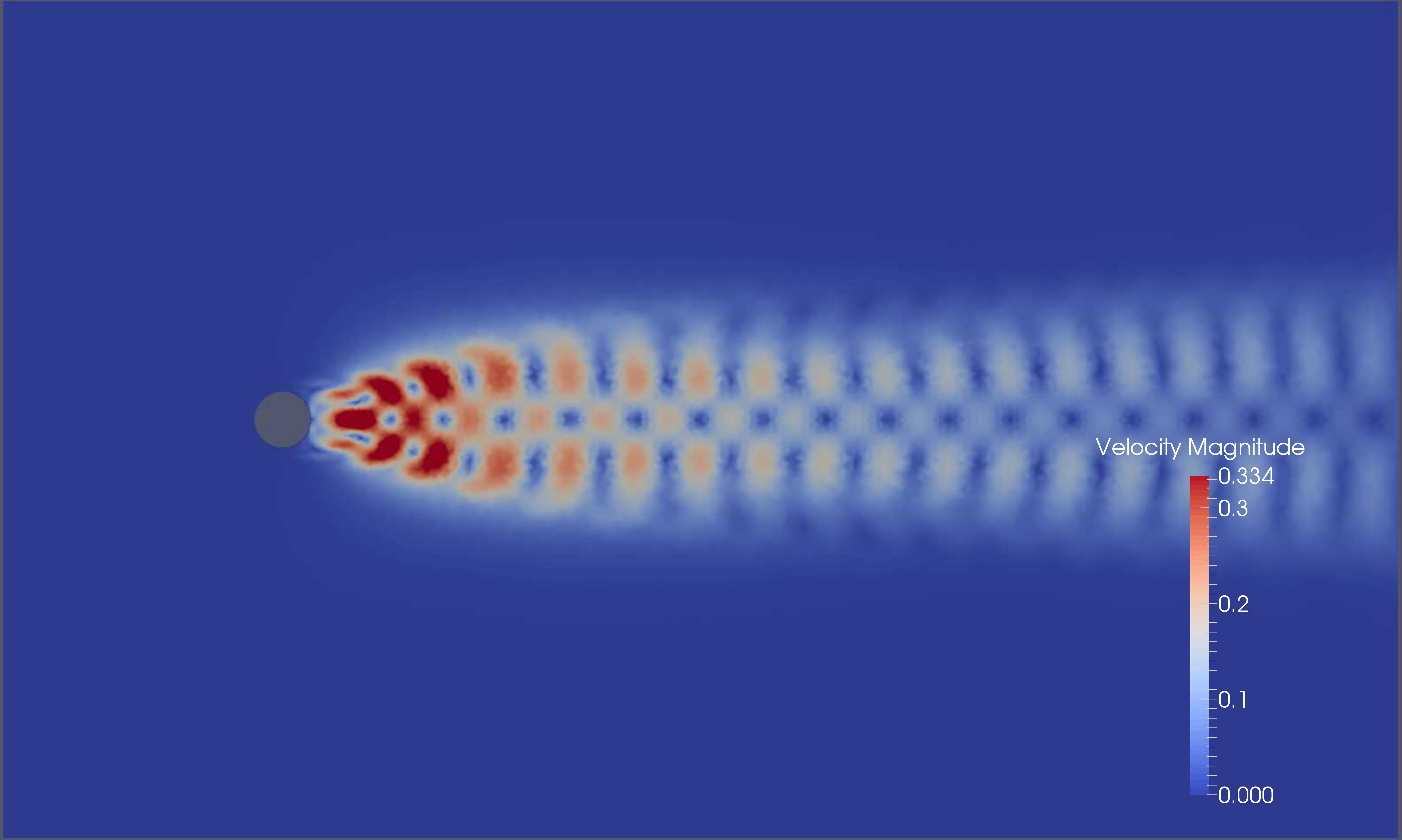}
	\caption{The first four POD modes for an uncontrolled solution.}
	\label{fig:POD_Modes}
\end{figure}

\begin{figure}[h!]
	\centering
	\parbox[b]{0.5\textwidth}{\centering \includegraphics[width=0.45\textwidth,height=0.25\textwidth]{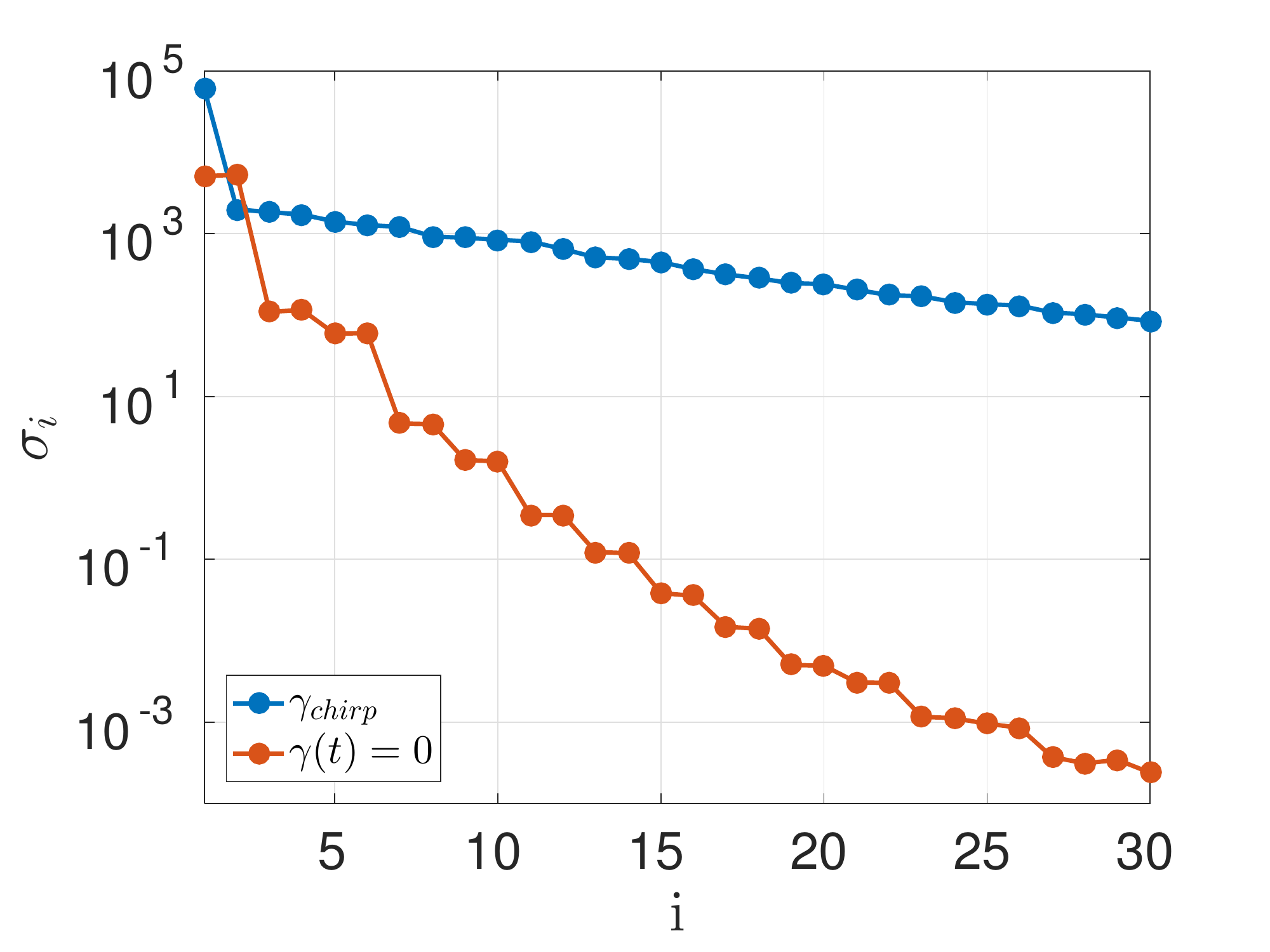}\\(a)}
	\parbox[b]{0.49\textwidth}{\centering \includegraphics[width=0.45\textwidth,height=0.25\textwidth]{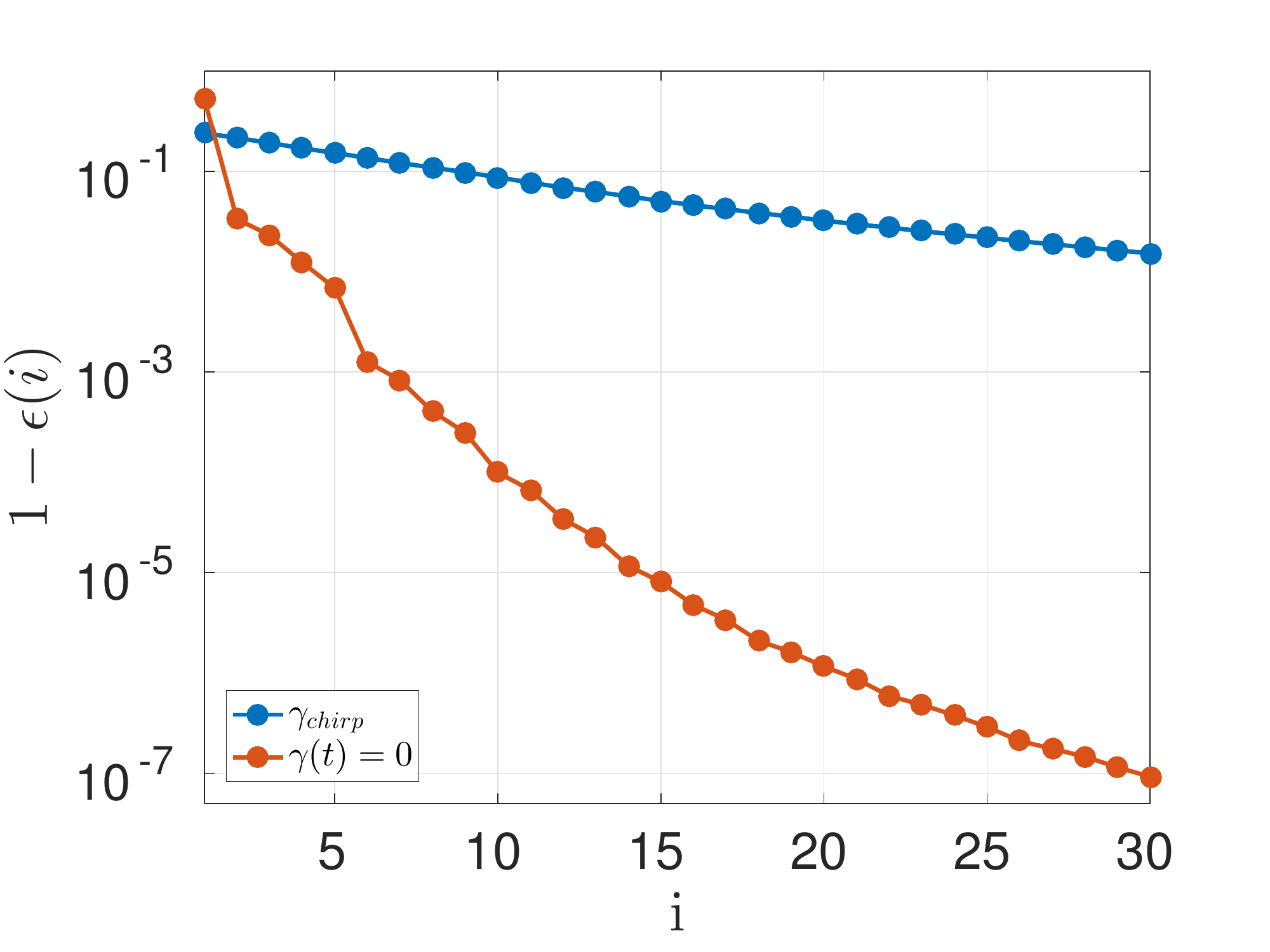}\\(b)}
	\caption{Eigenvalues (a) and error at current basis size (b) for an uncontrolled solution and a solution controlled by the chirping function.}
	\label{fig:Eigenvalues}
\end{figure}

Figure~\ref{fig:POD_Modes} shows the first four POD modes for the uncontrolled solution. The modes occur in pairs, the second one slightly shifted downstream. This is due to symmetries (see the double eigenvalues in Figure~\ref{fig:Eigenvalues}(a)) in the problem, namely in the horizontal axis through the cylinder as well as on the upper and lower boundary, respectively. The first four modes already account for $\approx 98 \%$ of the information, cf.~Figure~\ref{fig:Eigenvalues}, where in (a) the eigenvalues are shown for both the uncontrolled flow and the flow controlled by the chirping function. The error, i.e.~the ratio of the truncated eigenvalues, is visualized in (b).

Making use of the orthonormality of the POD basis, the first objective in \eqref{eq:R-MOCP1} can be simplified:
\begin{align*}
	\int_{\Omega} \Big \| \sum_{i=1}^l \alpha_i(t) \boldsymbol{\psi}_i(\boldsymbol{x}) \Big \|^2_2 \, d\boldsymbol{x} &= \int_{\Omega} \left( \sum_{i=1}^l \alpha_i(t) \boldsymbol{\psi}_i(\boldsymbol{x}) \right) \cdot \left(\sum_{i=1}^l \alpha_i(t) \boldsymbol{\psi}_i(\boldsymbol{x})\right) \, d\boldsymbol{x} = \sum_{i=1}^l \alpha_i^2(t) \\
	&\text{with } \left( \boldsymbol{\psi}_i, \boldsymbol{\psi}_j\right)_{L^2} = \delta_{i,j},
\end{align*}
which allows us to replace \eqref{eq:R-MOCP1} by a simpler formulation:
\begin{align}
	\min_{\gamma} J(\alpha, \gamma) = \min_{\gamma} \left( \begin{array}{c} \int_{t_0}^{t_e} \sum_{i=1}^l \alpha_i^2(t) \, dt \\ l\int_{t_0}^{t_e} \gamma^2(t) \, dt \end{array} \right).  \tag{R-MOCP} \label{eq:R-MOCP}
\end{align}

\subsection{Adjoint systems}
\label{subsec:Adjoint}
We would like to use a gradient based algorithm for the scalar optimization problems occurring in the reference point method (Section~\ref{subsubsec:MOCP_reference_point}). We realize this with an adjoint approach. The first step is to transform (\ref{eq:R-MOCP}) into a single objective optimization problem where the euclidean distance to a target $T$ is to be minimized:
\begin{align}
	\min_{\gamma} \overline{J}(\alpha, \gamma) &= \min_{\gamma} \|T - J(\alpha, \gamma) \|_2^2, \tag{SOP} \label{eq:SOP} \\
	\text{with} \hspace{1cm} \overline{J}(\alpha, \gamma) &= \left( T_1 - J_1(\alpha, \gamma) \right)^2 + \left( T_2 - J_2(\alpha, \gamma) \right)^2 \notag \\
	&= \left( T_1 - \int_{t_0}^{t_e} \sum_{i=1}^l \alpha_i^2(t) \, dt \right)^2 
	+ \left( T_2 - l\int_{t_0}^{t_e} \gamma^2(t) \, dt \right)^2, \label{eq:J_scalar}
\end{align}
where $J(\alpha,\gamma)$ according to \eqref{eq:R-MOCP} and $\alpha(t)$ satisfies the reduced order state equation \eqref{eq:eROM_state}, \eqref{eq:eROM_initial}. 

The task is now to solve a sequence of scalar optimization problems with varying targets $T$. We choose a line search strategy \cite{NW06} to address \eqref{eq:SOP}, where the direction is computed with a conjugate gradient method and the step length via backtracking and a valid length has to satisfy the Armijo rule.
The necessary gradient information is computed with an adjoint approach. To this end, we introduce the adjoint state $\lambda \in H^1([t_0,t_e],\mathbb{R}^l)$ and the Lagrange functional:
\begin{align}
	L(\alpha,\lambda, \gamma) = \int_{t_0}^{t_e} \overline{J}(\alpha, \gamma) - \lambda^\top \left( \dot{\alpha} - \mathcal{A} - \mathcal{B} \alpha - \mathcal{C}(\alpha) - \mathcal{D} \dot{\gamma} - (\mathcal{E} + \mathcal{F} \alpha) \gamma - \mathcal{G} \gamma^2 \right) \, dt, \label{eq:Lagrange1}
\end{align}
which is stationary for optimal values of $\gamma(t)$ and the corresponding state $\alpha(t)$ and adjoint state $\lambda(t)$. 
Using integration by parts, this leads to the following system of equations\footnote{Note, that since this approach makes use of variational calculus, it
is formally only applicable in the case of stronger assumptions, namely for $\alpha, \lambda, \gamma \in C^1$ (since we apply integration
by parts).}:
\begin{align}
	\dot{\alpha} &= \mathcal{A} + \mathcal{B} \alpha + \mathcal{C}(j\alpha) + \mathcal{D} \dot{\gamma} + (\mathcal{E} + \mathcal{F} \alpha) \gamma + \mathcal{G} \gamma^2, \label{eq:State1} \\
	\dot{\lambda} &= -\frac{\partial \overline{J}}{\partial \alpha} - \left(\mathcal{B} + \frac{\partial \mathcal{C}(\alpha)}{\partial \alpha} + \mathcal{F \gamma} \right)^\top \lambda, \label{eq:Adjoint1} \\
	0 &= \frac{\partial \overline{J}}{\partial \gamma} + \left(\mathcal{E} + \mathcal{F} \alpha + 2 \mathcal{G} \gamma \right)^\top \lambda - \mathcal{D}^\top \dot{\lambda} =: D_\gamma \overline{J}, \label{eq:Optimality1}
\end{align}
with the Fr\'echet derivatives
\begin{align*}
\partial \overline{J} / \partial \alpha_j &= 4 \left( \int_{t_0}^{t_e} \sum_{i=1}^l \alpha_i^2(t) \, dt - T_1 \right) \alpha_j(t), \\\partial \overline{J} / \partial \gamma &= 4 l \left( l \int_{t_0}^{t_e} \gamma^2(t) \, dt - T_2 \right) \gamma(t),
\end{align*}
and the respective boundary conditions $\alpha(0) = \alpha_0$ and $\lambda(t_e) = 0$. A detailed derivation of the optimality system is given in Appendix \ref{appendix:OS}. 

Equations \eqref{eq:State1} -- \eqref{eq:Optimality1} form the so-called optimality system which is seldom solved explicitly. Instead, the state equation \eqref{eq:State1} and the adjoint equation \eqref{eq:Adjoint1} are solved by forward and backward integration, respectively. The optimality condition \eqref{eq:Optimality1} provides the derivative information $D_\gamma \overline{J}$ of the cost functional with respect to the control $\gamma$. This information can be used to improve an initial guess for the optimal control within an iterative optimization scheme as described in Algorithm \ref{Alg:Backtracking}.
\begin{algorithm}
	\caption{(Adjoint based scalar optimization)}
	\label{Alg:Backtracking}
	\begin{algorithmic}[1]
	\Require Target point $T$, initial guess $\gamma^{(0)}$
	\For{$i=0,\ldots$}
	\State Compute $\alpha^{(i)}$ by solving the reduced order model \eqref{eq:State1} with the control $\gamma^{(i)}$
	\State Solve the adjoint equation \eqref{eq:Adjoint1} in backwards time with $\gamma^{(i)}$, $\alpha^{(i)}$, $\partial \overline{J} / \partial \alpha^{(i)}$
	\State Compute the gradient $D_\gamma \overline{J}^{(i)}$ by evaluating the optimality condition \eqref{eq:Optimality1} with $\gamma^{(i)}$, $\alpha^{(i)}$ and $\lambda^{(i)}$
	\State Update the control: $\gamma^{(i+1)} = \gamma^{(i)} + a^{(i)} d^{(i)}$, where $a^{(i)}$ is the step length computed by the backtracking Armijo method \cite{NW06} and $d^{(i)}$ is the conjugate gradient direction $d^{(i)}=-D_\gamma \overline{J}^{{(i)}} + \beta^i d^{{(i-1)}}$, $d^{(0)} = D_\gamma \overline{J}^{(0)}$ with $\beta^{(i)} = \frac{D_\gamma \overline{J}^{(i)} \, \cdot \, D_\gamma \overline{J}^{(i)}}{D_\gamma \overline{J}^{(i-1)} \, \cdot \, D_\gamma \overline{J}^{(i-1)}}$
	\EndFor
	\end{algorithmic}
\end{algorithm}

\begin{figure}[h!]
	\centering
	\parbox[b]{0.5\textwidth}{\centering \includegraphics[width=0.5\textwidth]{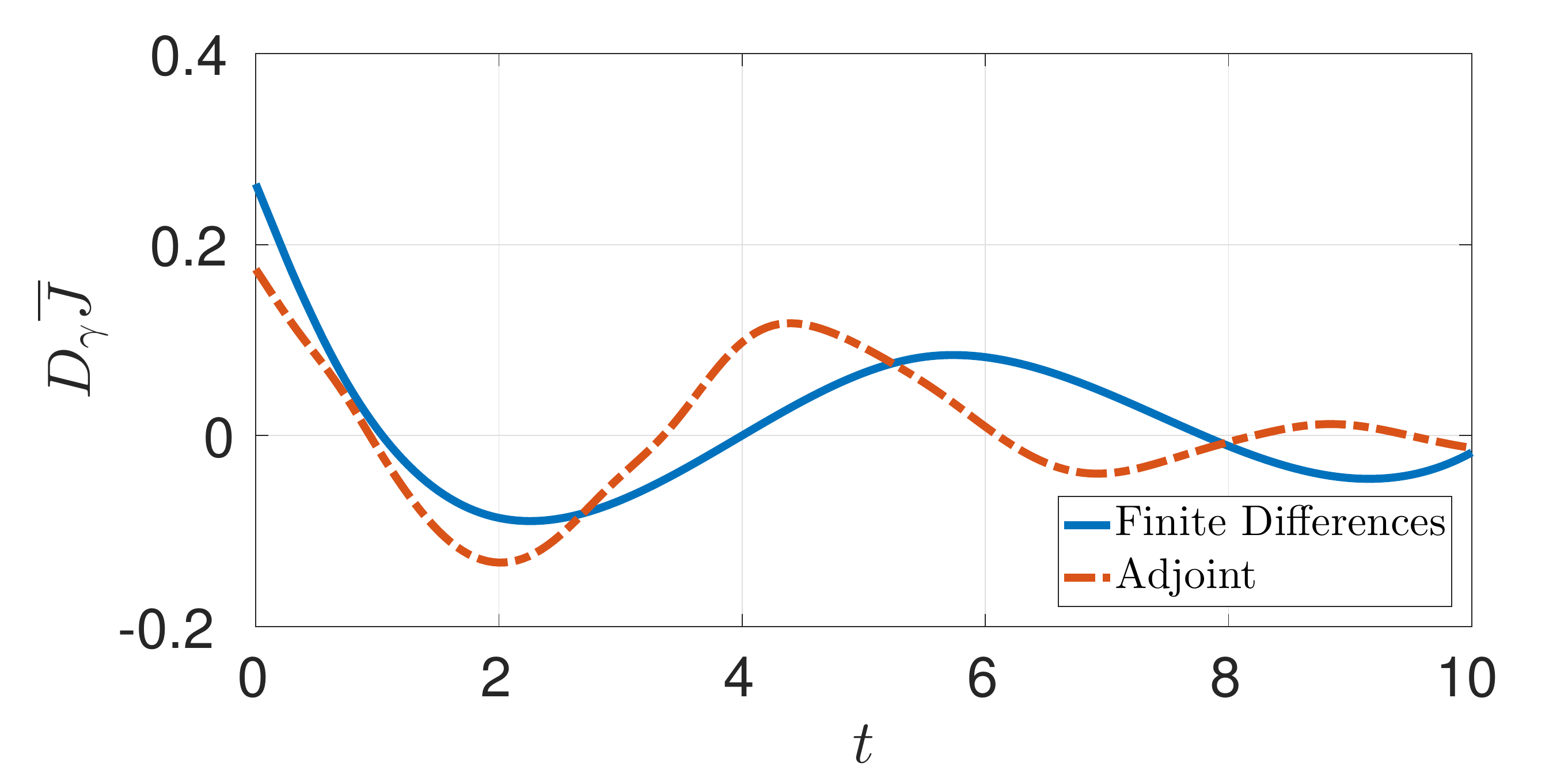}\\(a)}
	\parbox[b]{0.49\textwidth}{\centering \includegraphics[width=0.49\textwidth]{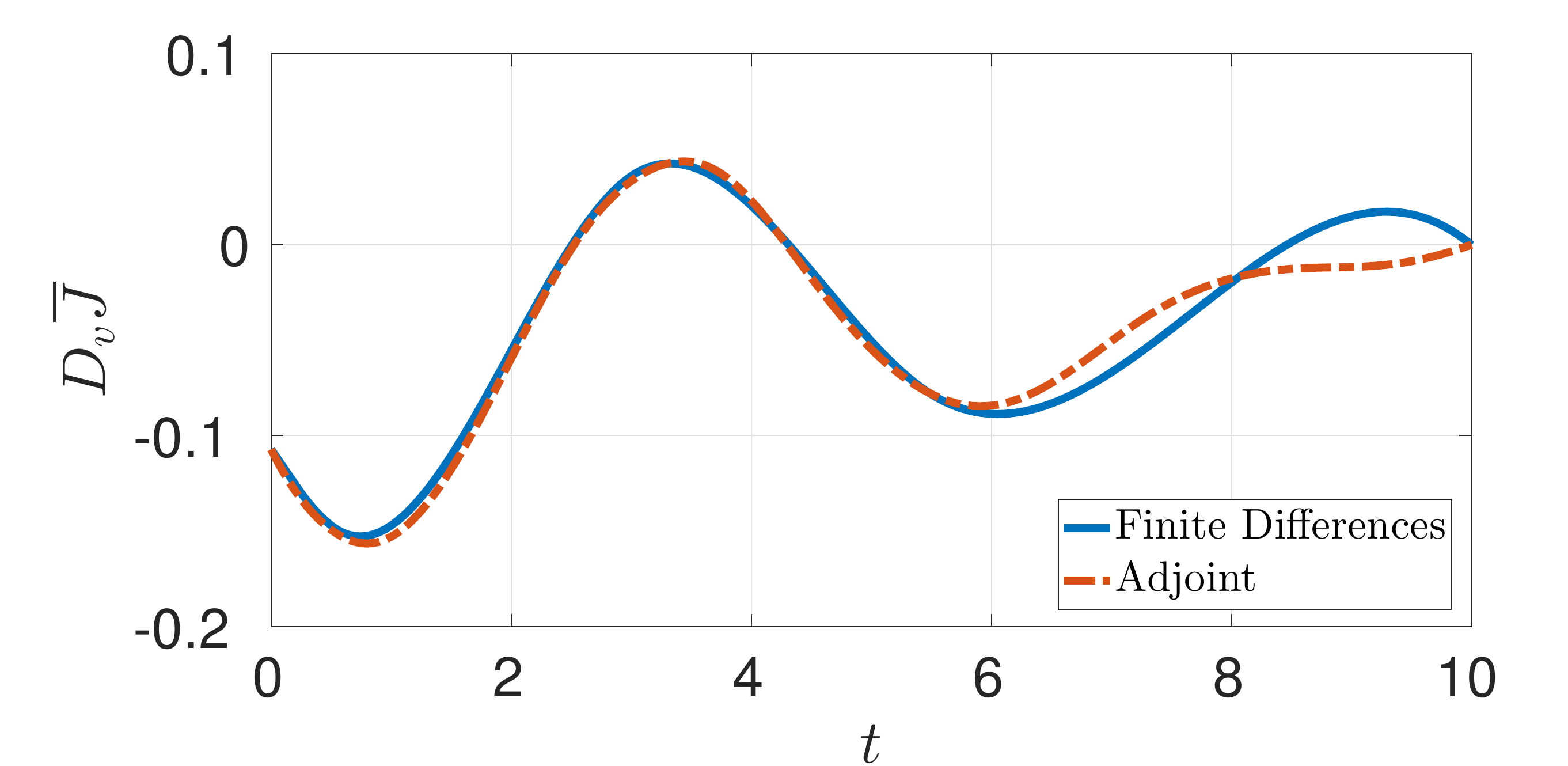}\\(b)}
	\caption{Comparison of the gradient computed using the adjoint approaches with an approximation by finite differences (FD). (a) Optimality system \eqref{eq:State1} -- \eqref{eq:Optimality1} (FD: $D_{\gamma,i} \overline{J} = (\overline{J}(\gamma + h e_i) - \overline{J}(\gamma))/h$, where $h e_i$ is the variation by $h$ of the $i$-th entry of the discretization of $\gamma$). (b) Optimality system \eqref{eq:State2_a} -- \eqref{eq:Optimality2} (FD: $D_{v,i} \overline{J} = (\overline{J}(v + h e_i) - \overline{J}(v))/h$).}
	\label{fig:Gradients_FD_vs_Adj}
\end{figure}

\begin{figure}
	\centering
	\includegraphics[width=0.4\textwidth]{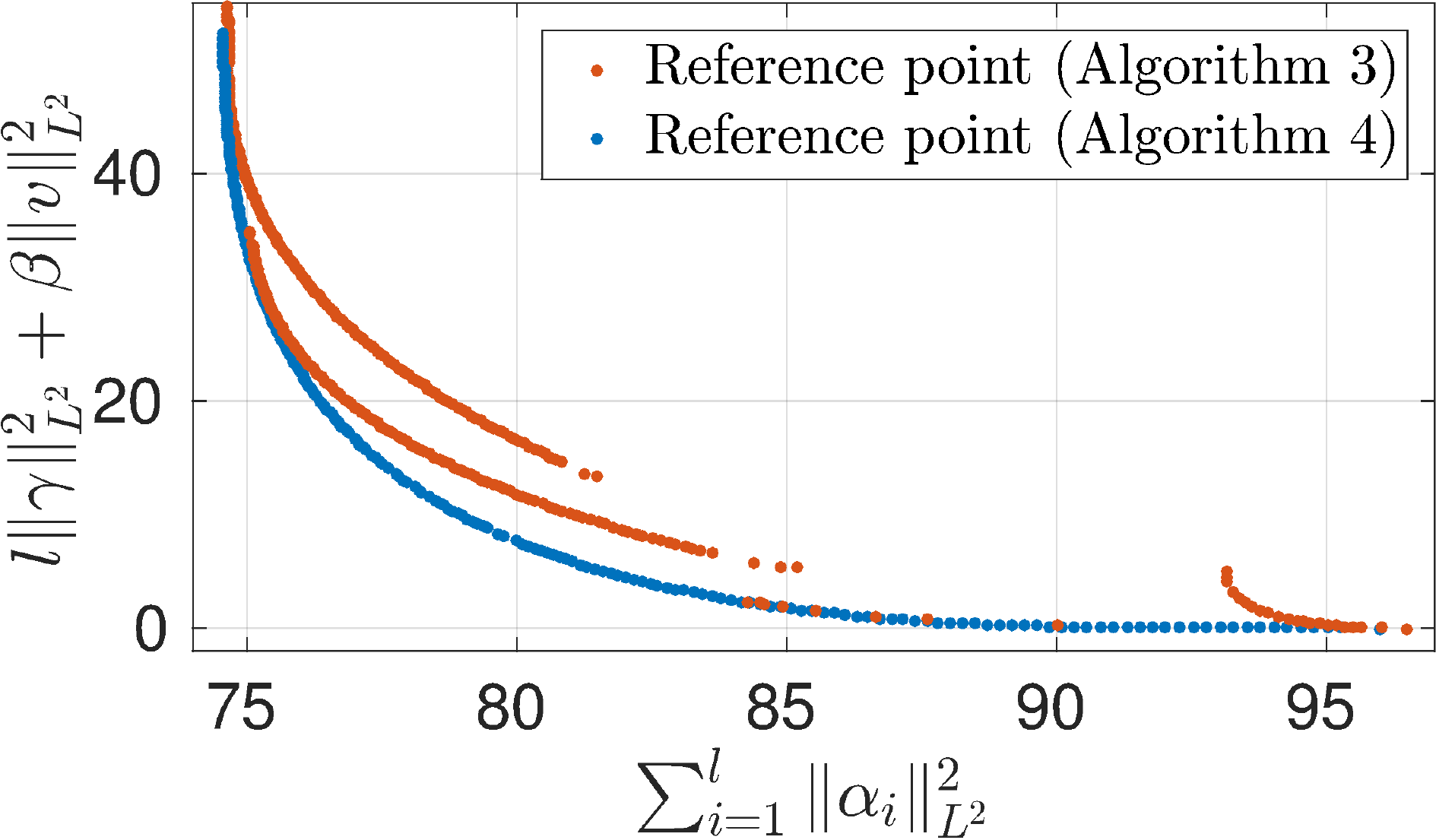}
	\caption{Comparison of the Pareto fronts computed with the two optimality systems \eqref{eq:State1} -- \eqref{eq:Optimality1} and \eqref{eq:State2_a} -- \eqref{eq:Optimality2}.}
	\label{fig:Comparison_optimality_systems}
\end{figure}

In \cite{BCB05}, the condition $\lambda(t_e) = 0$ is the only boundary condition for the adjoint equation and the condition $\lambda^\top(t_0) \mathcal{D} \delta \gamma(t_0) = 0$ (cf.~Appendix~\ref{appendix:OS}) is neglected. In our computations, using this optimality system caused convergence issues. The reason for this is that the gradient obtained by the adjoint approach appears to be wrong, cf.~Figure~\ref{fig:Gradients_FD_vs_Adj}(a) for a comparison to a finite difference approximation. (Since the state equation is very accurate, the gradient can be approximated well this way.) The reason for this might be that we neglect the additional condition for the adjoint equation. However, it has been reported before \cite{Fah00, DV01} that it is favorable to include information of the adjoint state of the PDE in the reduced order model. In fact, the adjoint based gradient can easily become inaccurate or even incorrect if based on information about the state only. Consequently, it is difficult to say whether the convergence problems stem from the missing boundary condition or are a flaw of the model reduction approach itself. To further investigate this, we derive a second optimality system based on an alternative formulation of the reduced model. Following \cite{Rav00}, we define an augmented state $(\alpha, \gamma)^\top$ and introduce a new control $v \in L^2([t_0,t_e], \mathbb{R}),\ v(t) = \dot{\gamma}(t)$. 
It is known from optimal control theory that if a dynamical system depends linearly on the control, the control is bounded by box constraints and does not appear in the cost functional, the optimal control is of \emph{bang bang} type and might include singular arcs \cite{Ger12}. This is exactly the situation for the second formulation and since the computation of such solutions is numerically challenging, we avoid this behavior by adding a regularization term to the second objective in \eqref{eq:R-MOCP}, i.e.~$J_2 = l \|\gamma\|_{L^2}^2 + \beta \|v\|_{L^2}^2$, where $\beta$ is a small number (here: $\beta = 1e^{-5}$). The optimality system can be derived in an analogous way to the one described above, where the adjoint states are now denoted by $\lambda \in H^1([t_0,t_e],\mathbb{R}^l)$ and $\mu \in H^1([t_0,t_e],\mathbb{R})$:
\begin{align}
 	\dot{\alpha} &= \mathcal{A} + \mathcal{B} \alpha + \mathcal{C}(\alpha) + \mathcal{D} v + (\mathcal{E} + \mathcal{F} \alpha) \gamma + \mathcal{G} \gamma^2, \label{eq:State2_a} \\
 	\dot{\gamma} &= v, \label{eq:State2_b} \\
	\dot{\lambda} &= -\frac{\partial \overline{J}}{\partial \alpha} - \left(\mathcal{B} + \frac{\partial \mathcal{C}(\alpha)}{\partial \alpha} + \mathcal{F \gamma} \right)^\top \lambda, \label{eq:Adjoint2_a} \\
	\dot{\mu} &= -\frac{\partial \overline{J}}{\partial \gamma} - \left(\mathcal{E} + \mathcal{F \alpha} + 2 \mathcal{G} \gamma \right)^\top \lambda, \label{eq:Adjoint2_b} \\
	0 &= \frac{\partial \overline{J}}{\partial v} + \mathcal{D}^\top \lambda + \mu =: D_v \overline{J}, \label{eq:Optimality2}
\end{align}
with the Fr\'echet derivatives:
\begin{align*}
\partial \overline{J} / \partial \gamma &= 4 l \left( l \int_{t_0}^{t_e} \gamma^2(t) \, dt + \beta \int_{t_0}^{t_e} v^2(t) \, dt - T_2 \right) \gamma(t), \\
\partial \overline{J} / \partial v &= 4 \beta \left( l \int_{t_0}^{t_e} \gamma^2(t) \, dt + \beta \int_{t_0}^{t_e} v^2(t) \, dt - T_2 \right) v(t),
\end{align*} 
$\partial \overline{J} / \partial \alpha_j$ as before and the respective boundary conditions
\begin{align}
	\alpha(0) = \alpha_0, \ \lambda(T) = 0, \ \mu(0) = 0, \ \mu(T) = 0. \label{eq:IC_BC2}
\end{align}
This yields a boundary value problem for the state equations \eqref{eq:State2_a}, \eqref{eq:State2_b} and the adjoint equations \eqref{eq:Adjoint2_a}, \eqref{eq:Adjoint2_b} which can, for a given control $v$, be solved by an iterative scheme. Applying a shooting method, the initial value $\gamma(0)$ is computed such that $\mu(0) = 0$ is satisfied using Newton's method (cf.~Algorithm \ref{Alg:Backtracking_Shooting}). This is computationally more expensive than the simple forward-backward integration in Algorithm \ref{Alg:Backtracking}. However, the system \eqref{eq:State2_a} -- \eqref{eq:Optimality2} yields strongly improved convergence and decreased sensitivity to the optimization parameters in comparison to the system \eqref{eq:State1} -- \eqref{eq:Optimality1}, which stems form the improved gradient accuracy (cf.~Figure~\ref{fig:Gradients_FD_vs_Adj}(b)). This is also evident from Figure~\ref{fig:Comparison_optimality_systems}, where the Pareto front based on \eqref{eq:State2_a} -- \eqref{eq:Optimality2} was computed executing Algorithm~\ref{Alg:Backtracking_Shooting} once, starting in the already known point of zero control cost. In contrast to that, multiple attempts were made with Algorithm~\ref{Alg:Backtracking}, starting from different points on the Pareto front that were computed with Algorithm~\ref{Alg:Backtracking_Shooting}. All of them either divert quickly from the front or stop prematurely. The results shown in Section~\ref{sec:Results} are therefore based on \eqref{eq:State2_a} -- \eqref{eq:Optimality2}. From our experience, the approach based on the adjoint equation can be very sensitive, also with respect to stabilization methods where the coefficients are manipulated to better match the state equation. Hence, it would be favorable for future work to utilize direct approaches where only accuracy of the state equation is of importance.

\begin{algorithm}
	\caption{(Adjoint based scalar optimization with shooting)}
	\label{Alg:Backtracking_Shooting}
	\begin{algorithmic}[1]
	\Require Target point $T$, initial guess $v^{0}$
	\For{$i=0,\ldots$}
	\State Shooting step: Determine $\gamma^{(i)}(0)$ by solving an internal root finding problem in order to enforce the condition $\mu^{(i)}(0) = 0$. This step requires forward solves of \eqref{eq:State2_a}, \eqref{eq:State2_b} and backward solves of \eqref{eq:Adjoint2_a}, \eqref{eq:Adjoint2_b}.
	\State Compute $\alpha^{(i)}$, $\gamma^{(i)}$ by solving the reduced order model \eqref{eq:State2_a}, \eqref{eq:State2_b} with $v^{(i)}$ as the control and $\gamma^{(i)}(0)$ as computed in step 2.
	\State Solve the adjoint equations \eqref{eq:Adjoint2_a}, \eqref{eq:Adjoint2_b} in backwards time with $v^{(i)}$, $\alpha^{(i)}, \gamma^{(i)}$, $\left(\partial \overline{J} / \partial \alpha\right)^{(i)}$ and $\left(\partial \overline{J} / \partial \gamma\right)^{(i)}$.
	\State Compute the gradient $\left(D_v \overline{J}\right)^{(i)}$ by evaluating the optimality condition \eqref{eq:Optimality2} with $\lambda^{(i)}$  and $\mu^{(i)}$.
	\State Update the control: $v^{(i+1)} = v^{(i)} + a^{(i)} d^{(i)}$, with $a^{(i)}$ and $d^{(i)}$ as in Algorithm \ref{Alg:Backtracking}.
	\EndFor
	\end{algorithmic}
\end{algorithm}

%% file: 05_Results.tex
In this section we present solutions to \eqref{eq:R-MOCP} obtained with the algorithms presented in Section~\ref{sec:MOC}. For all cases, we fix the time interval $[t_0, t_e]$ to $[0,10]$ which corresponds to roughly two vortex shedding cycles. In order to use the subdivision method, we transform \eqref{eq:MOCP} into a multiobjective optimization problem with moderate decision space dimension ($m \leq 10$). One way to do this is by introducing a sinusoidal control (cf.~Section~\ref{subsubsec:MOCP_subdivision}) which then yields \eqref{eq:MOP-2D}. Another way is to approximate $\gamma$ by a cubic spline with $m$ break points. In this case, the parameters to be optimized are the spline values at these break points. For the reference point method, we discretize $\gamma$ with a constant step length of $\Delta t = 0.05$ and solve Equations \eqref{eq:State2_a} -- \eqref{eq:Adjoint2_b} with a fourth order Runge-Kutta scheme.

\begin{figure}[h!]
	\includegraphics[width=0.5\textwidth,height=0.28\textwidth]{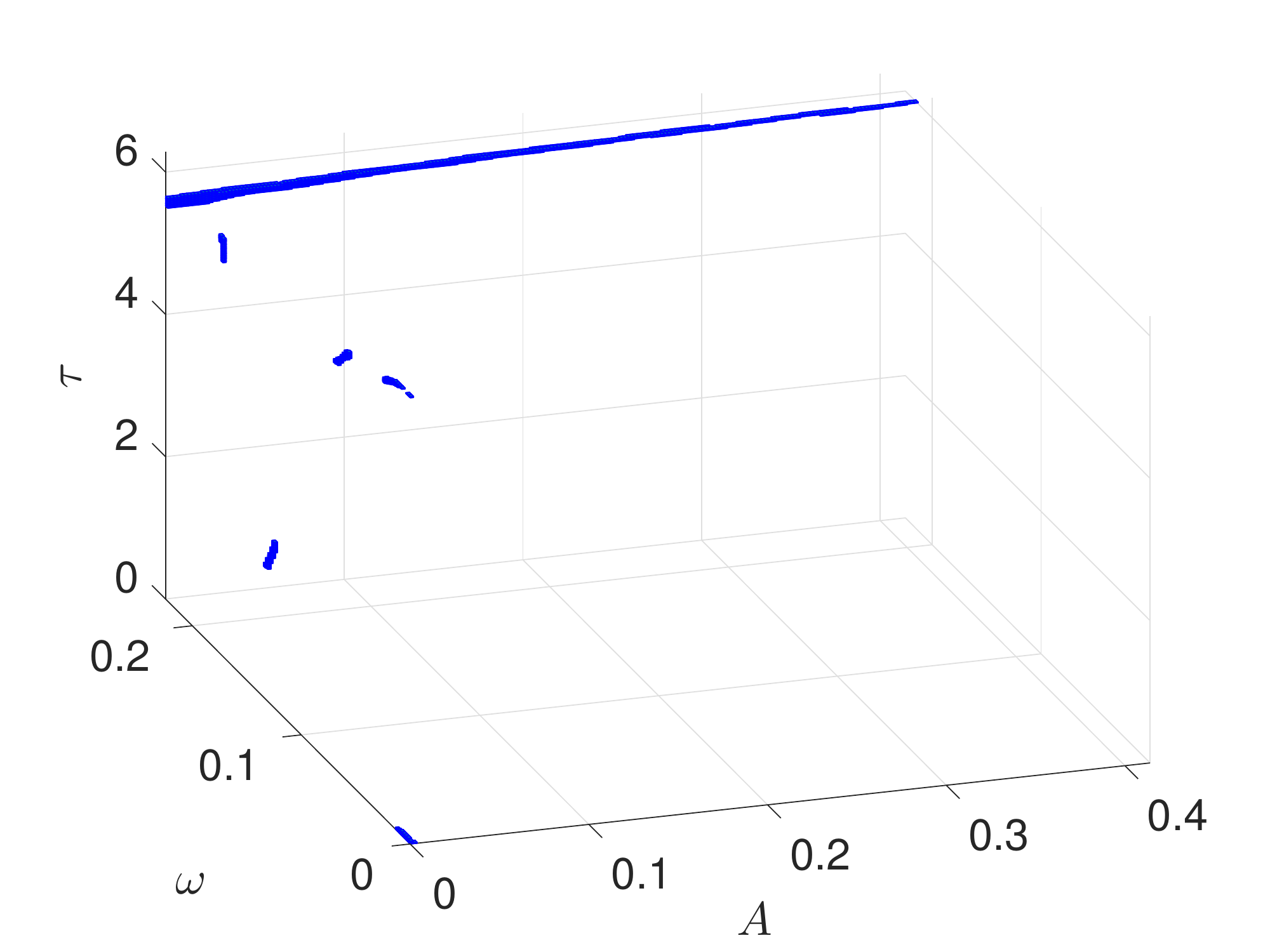}
	\includegraphics[width=0.5\textwidth,height=0.28\textwidth]{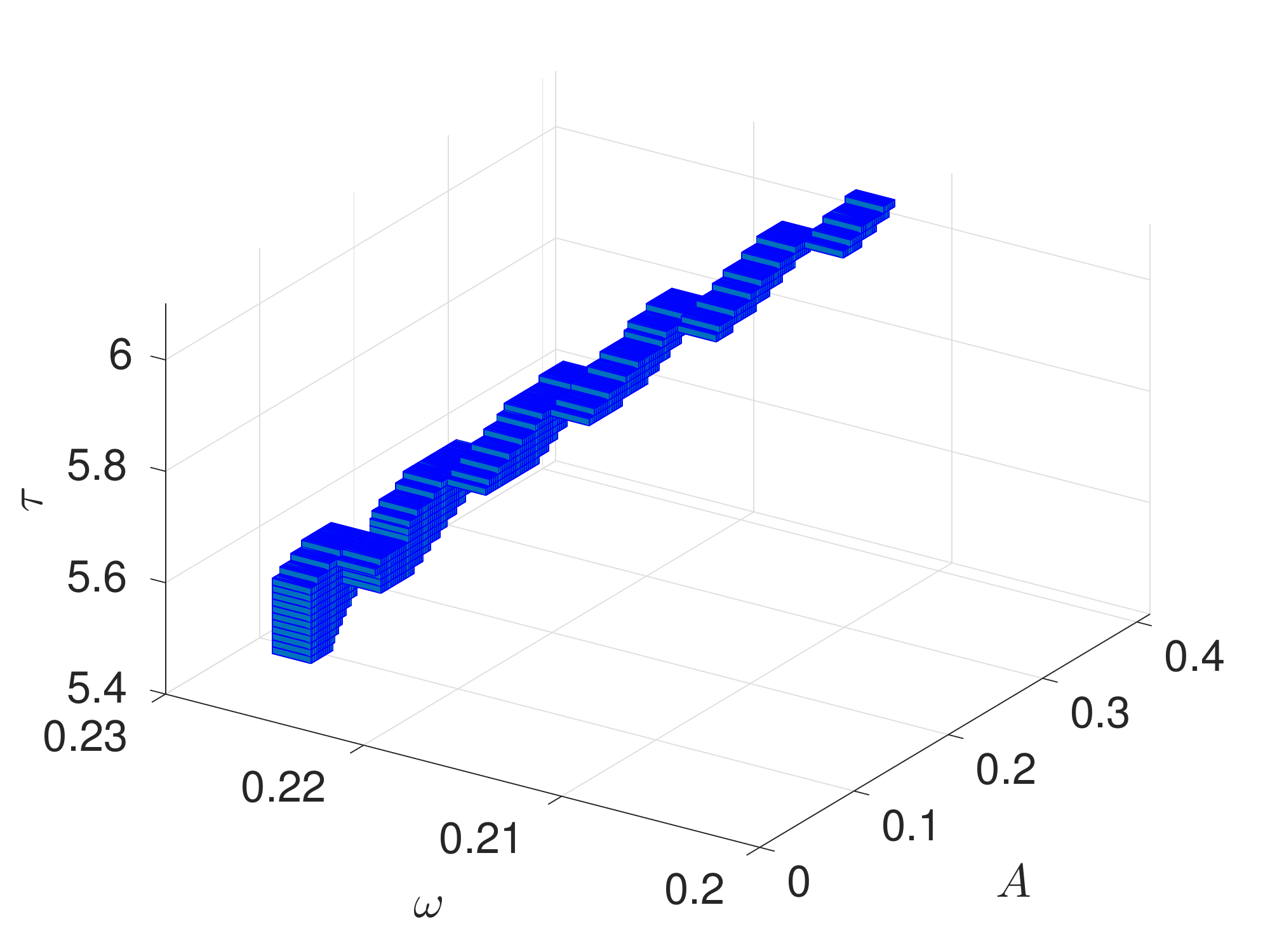}
	\caption{Box covering of the solution of \eqref{eq:MOP-2D}, obtained with global subdivision algorithm after 27 subdivision steps. There are spurious solutions with almost zero amplitude, the right plot shows only the physically relevant part which is constrained to a small part of the frequency domain. The trade-off between the two objectives is mainly realized by changing the amplitude.}
	\label{fig:results_PS_GAIO2D}
\end{figure}

\begin{figure}[t]
	\centering
	\parbox[b]{0.50\textwidth}{\centering \includegraphics[width=0.48\textwidth,height=0.28\textwidth]{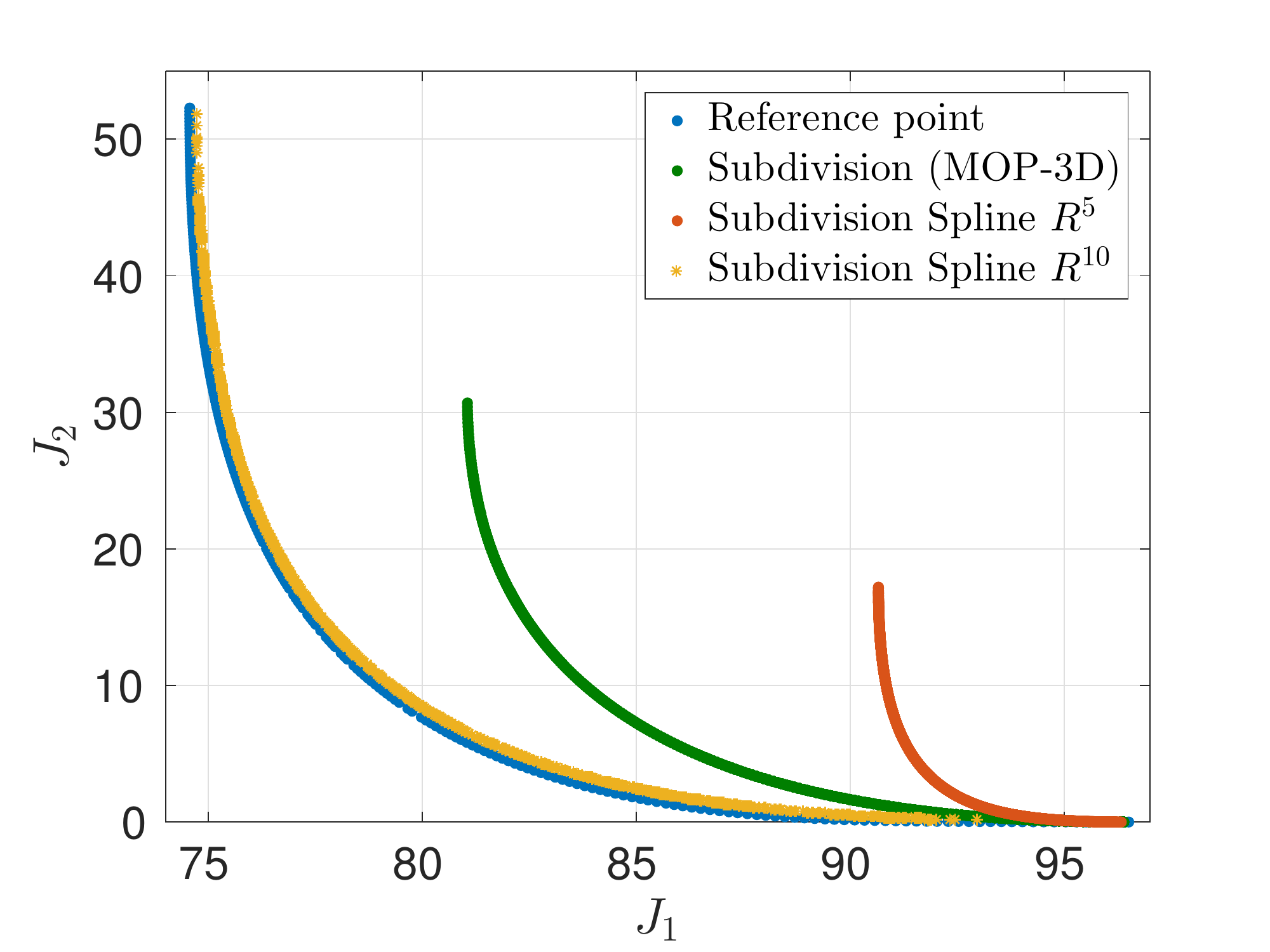}\\(a)}
	\parbox[b]{0.40\textwidth}{\centering \includegraphics[width=0.40\textwidth,height=0.28\textwidth]{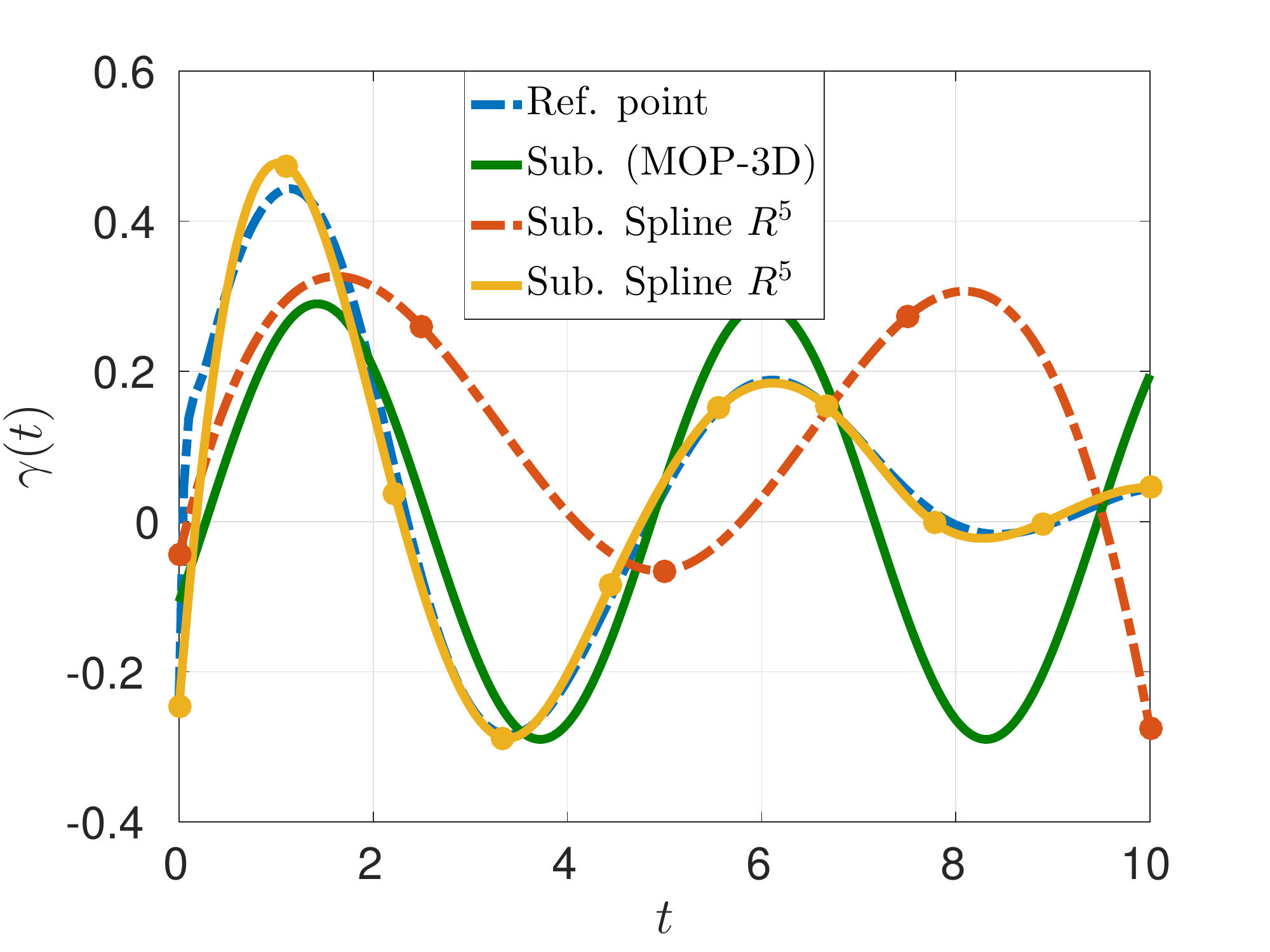}\\(b)}
	\caption{(a) Comparison of Pareto fronts for different MOCP solution methods. (b) Comparison of the Pareto optimal controls with equal cost ($J_2 = 15$) for different MOCP solution methods.}
	\label{fig:results_ParetoFronts}
\end{figure}

\begin{figure}[t]
	\centering
	\parbox[b]{0.49\textwidth}{\centering \includegraphics[width=0.45\textwidth,height=0.27\textwidth]{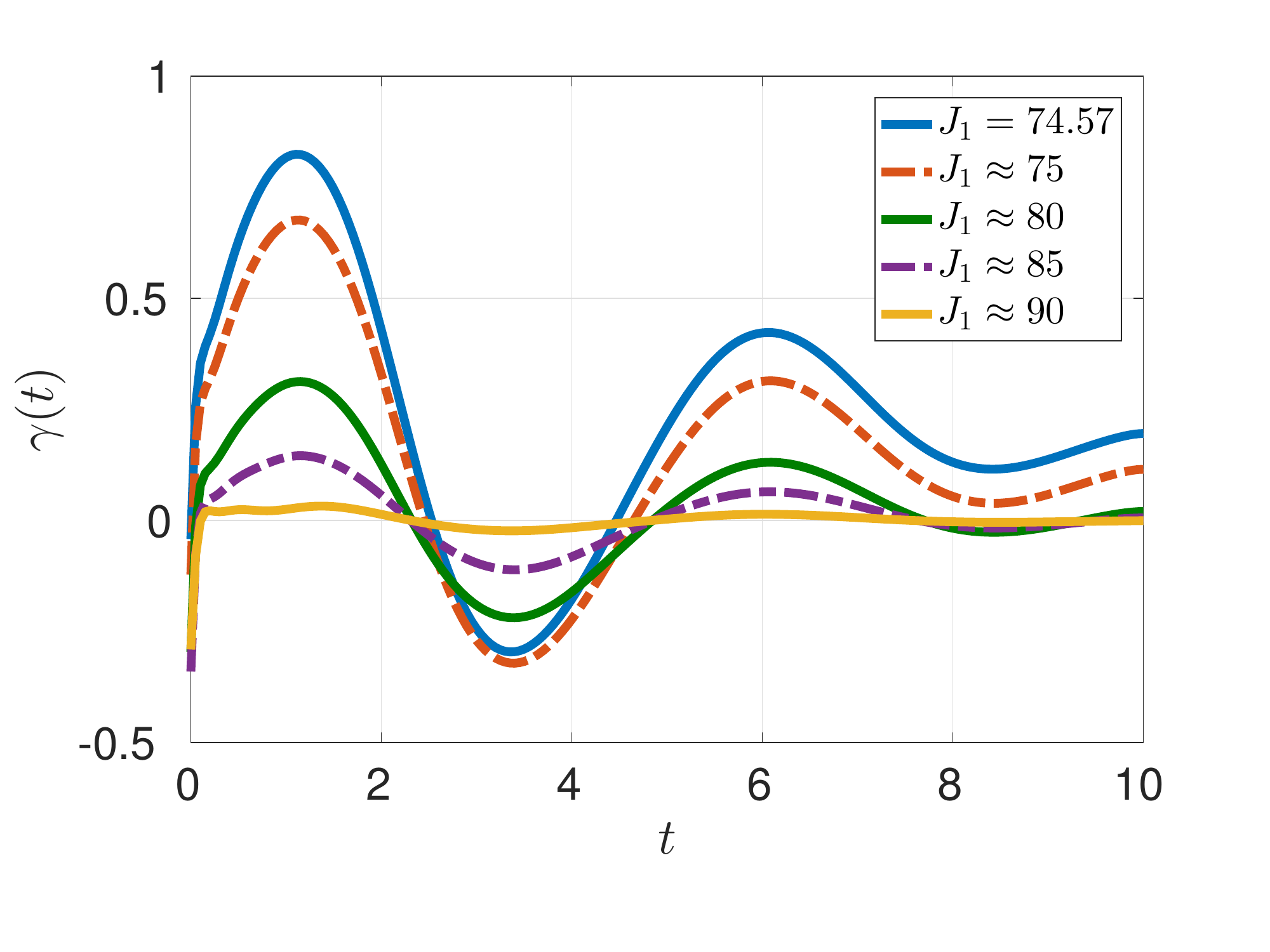}\\(a)}
	\parbox[b]{0.49\textwidth}{\centering \includegraphics[width=0.45\textwidth,height=0.27\textwidth]{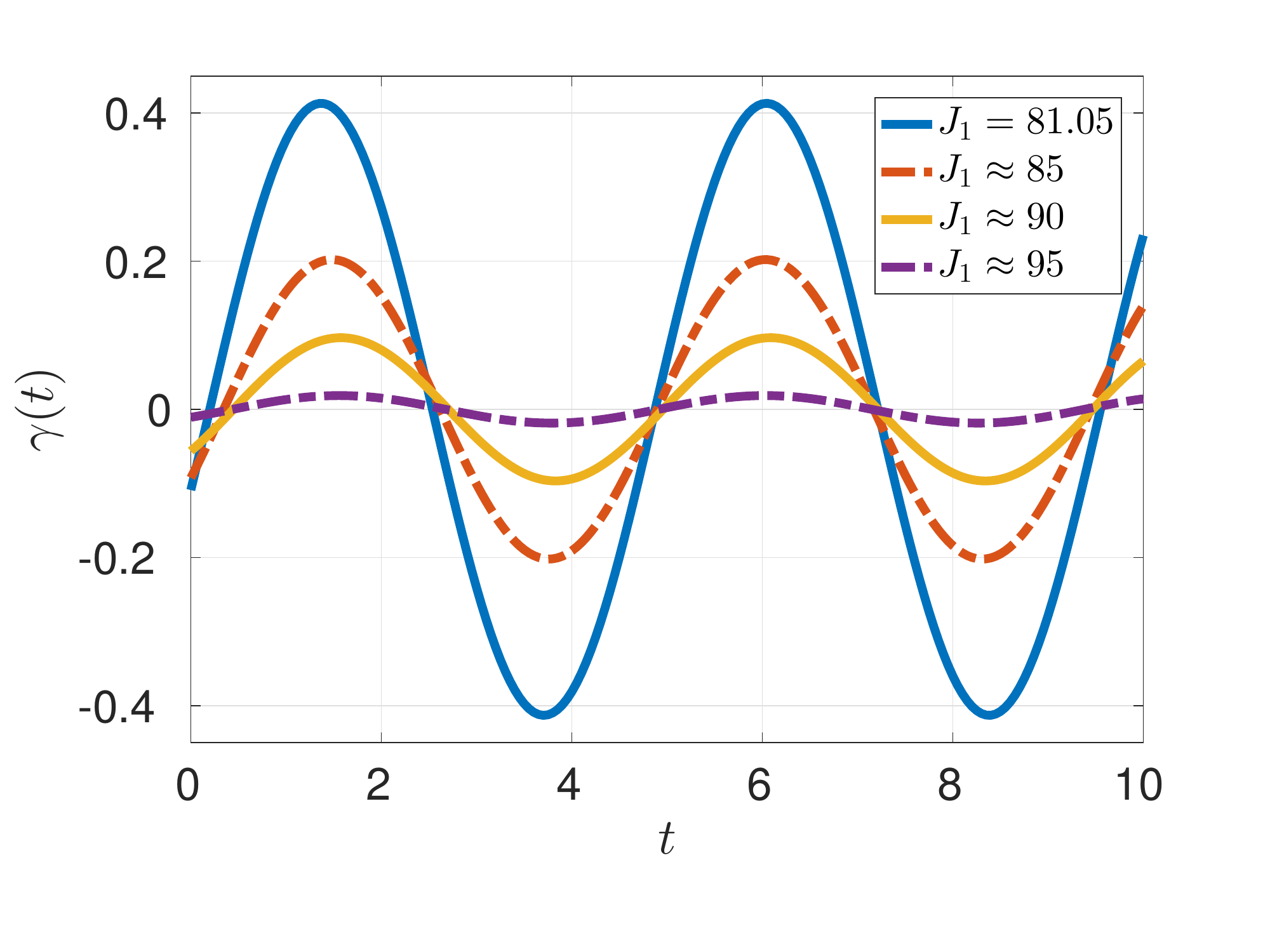}\\(b)}
	\caption{Different Pareto points for the reference point method (a) and \eqref{eq:MOP-2D} (b).}
	\label{fig:results_ParetoPoints_RPmethod}
\end{figure}

Figure~\ref{fig:results_PS_GAIO2D} shows the box covering of the Pareto set for \eqref{eq:MOP-2D}, the respective Pareto front is shown in Figure~\ref{fig:results_ParetoFronts}. It is noteworthy that, apart from a few spurious points at $A\approx 0$ caused by numerical errors (Figure~\ref{fig:MOCP_GAIO}(a)), the largest part of the set is restricted to a small part of the frequency domain. This is comprehensible since the rotation counteracts the natural dynamics of the von K\'{a}rm\'{a}n vortex street. The trade-off between control cost and stabilization is almost exclusively realized by changing the amplitude of the rotation.
When moving towards the two scalar optimal solutions, i.e.~minimal fluctuations and minimal control cost, respectively, further improvements are only possible by accepting large trade-offs in the other objective which is a common phenomenon in multiobjective optimization. When designing a system, one could now accept a small increase in the main objective, i.e.~flow stabilization, in order to achieve a large decrease in the control cost. This becomes more clear when looking at Figure~\ref{fig:results_ParetoPoints_RPmethod}(b), where different optimal compromises for the control are shown. For a relatively small improvement  of $J_1$ from $85$ to $81.05$ ($\approx -5\%$), the control cost increases by $322\%$.

\begin{figure}[t]
	\centering
	\parbox[b]{0.49\textwidth}{\centering \includegraphics[width=0.45\textwidth]{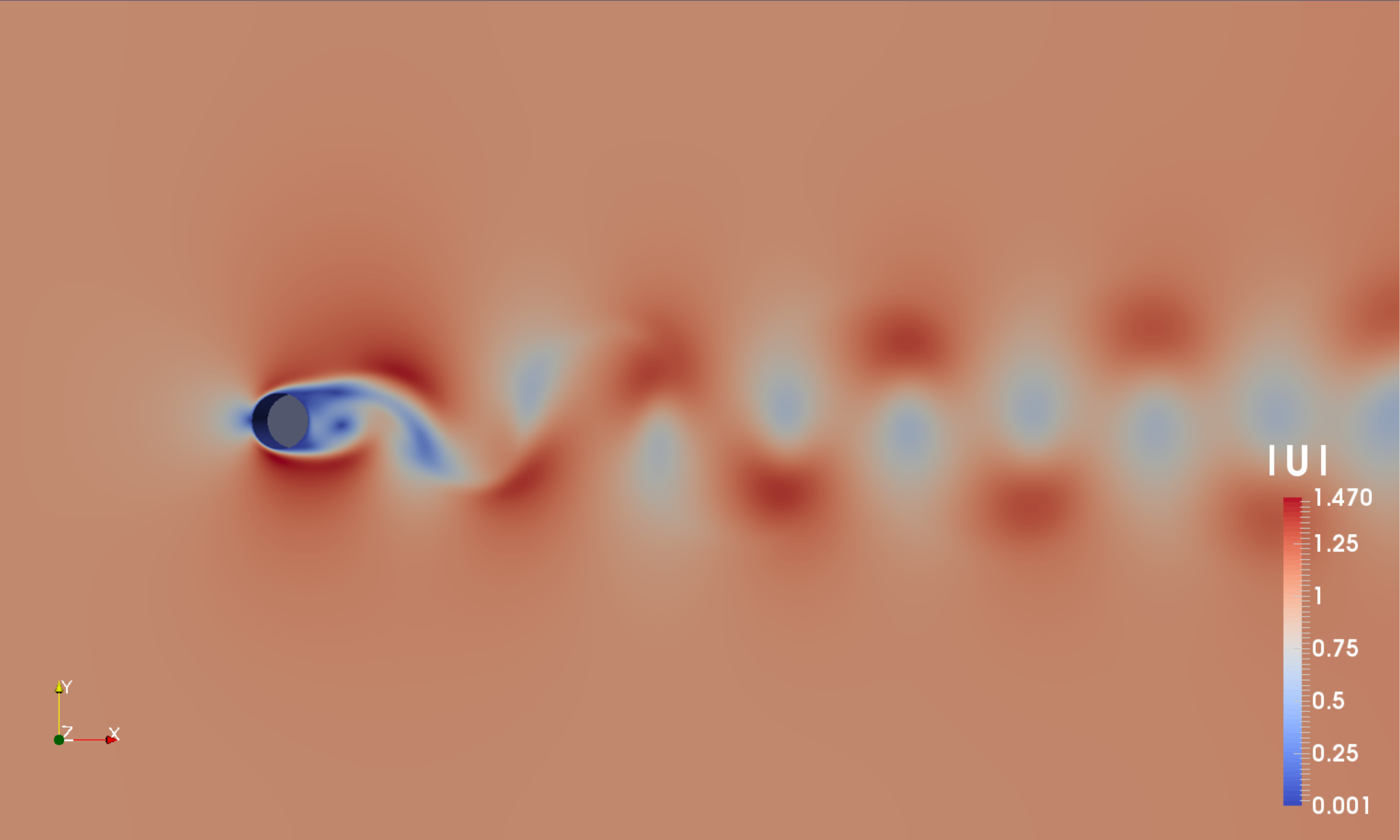}\\(a)}
	\parbox[b]{0.49\textwidth}{\centering \includegraphics[width=0.45\textwidth]{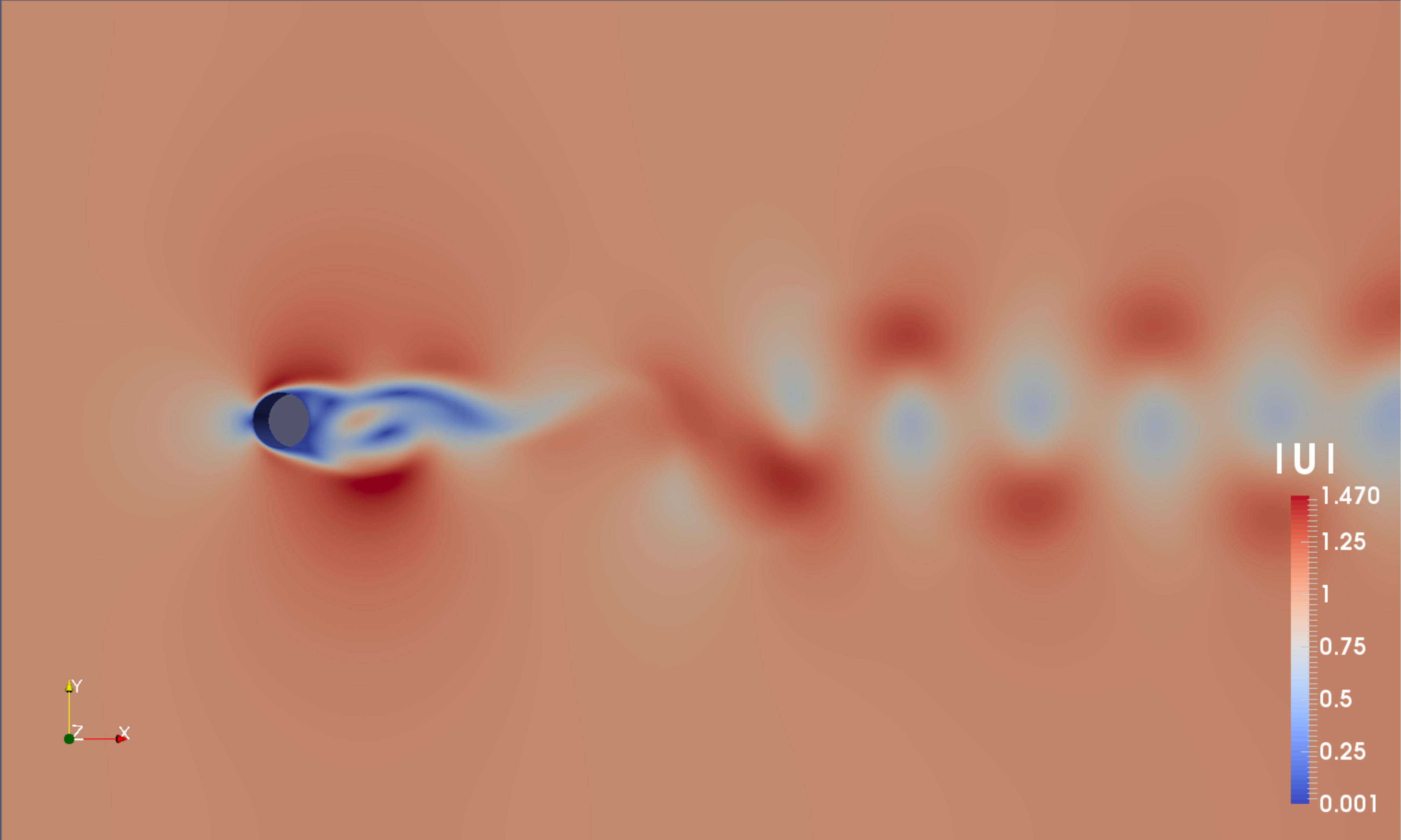}\\(b)}
	\caption{FV solution at $t = 10$ with two different Pareto optimal solutions computed with the reference point method. (a) Relatively low control cost ($J_2 = 1.84$), almost no reduction of fluctuations ($J_1=85$). (b) Larger control cost ($J_2 = 30.12$), stronger reduction of fluctuations ($J_1=75$). Due to the short integration time, the effect is only visible in the vicinity of the cylinder.}
	\label{fig:results_Solutions_RPmethod}
\end{figure}

Figures~\ref{fig:results_ParetoFronts}(a) and \ref{fig:results_ParetoPoints_RPmethod} also show the Pareto fronts and Pareto optimal controls, respectively, computed with the other methods, i.e.~two different spline approximations ($m = 5$ and $m = 10$) computed with the subdivision algorithm on the one hand and arbitrary control functions determined by the reference point method on the other hand. We observe that a spline with 5 breakpoints is too restrictive to find control functions of acceptable quality. When using 10 breakpoints on the other hand, the Pareto front clearly surpasses the solution of \eqref{eq:MOP-2D}. It is, however, numerically expensive to compute. The best solution is computed with the reference point method. The improvement over the spline based solution is relatively low, the reason being that the controls are almost similar (see Figure~\ref{fig:results_ParetoFronts}(b), where solutions of the different methods with the same control cost $J_2 = 15$ are compared). When considering scenarios with larger time intervals than 10 seconds, the spline dimension would have to be increased further, leading to again much higher computational cost whereas the cost for the reference point method increases only linearly with time due to the fact that only the number of time steps in the forward and backward integration \eqref{eq:State2_a} -- \eqref{eq:Adjoint2_b} is affected. Finally, Figure~\ref{fig:results_Solutions_RPmethod} shows two solutions of \eqref{eq:NSE1} -- \eqref{eq:BCC} with controls computed by the reference point method, corresponding to different values of the control cost $J_2$ and hence different degrees of stabilization. Due to the relatively short integration time of $t_e=10$ seconds, changes are only apparent closely past the cylinder. Still, we observe increased stabilization when allowing higher control cost.

\begin{figure}[h!]
	\centering
	\includegraphics[width=0.5\textwidth]{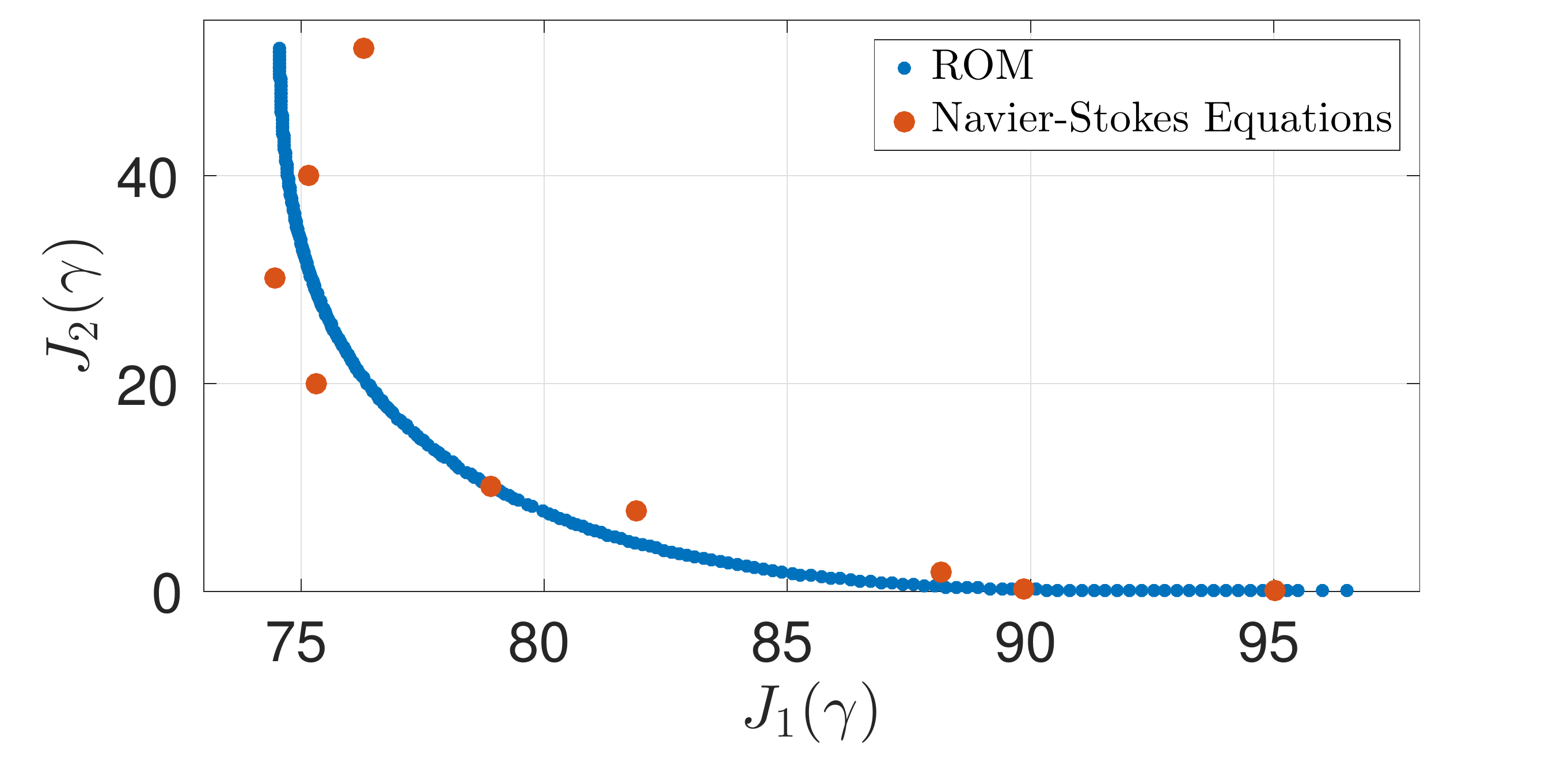}
	\caption{Comparison of the objective values based on the state of the reduced model and of the PDE solution. Since there is no error in $J_2$, the points are only shifted horizontally. Although the overall agreement is acceptable, the PDE based Pareto front appears to have a different shape, especially for solutions with larger control cost.}
	\label{fig:Validation_PDE_vs_ROM}
\end{figure}

In order to evaluate the quality of the solution obtained with the reduced order modeling approach, we have evaluated the cost function using the system state obtained from a PDE evaluation. This is depicted in Figure~\ref{fig:Validation_PDE_vs_ROM}. Since we only introduce an error in the first objective by the reduced order model, the error results in a horizontal shift of the Pareto points. In general, we observe a good agreement between the solutions, the error being less than $4\%$ for all points that were tested. However, especially in regions with a steep gradient $\partial J_2 / \partial J_1$ in the Pareto front, we see that this small error can result in a different shape of the front. Note that the two points with the highest values for $J_2$ are now dominated by the next lower point. This gives a strong motivation for further efforts to decrease the error of the reduced order model (e.g.~by Trust-Region methods \cite{Fah00}) and investigate the influence of the inaccuracies in the gradient obtained by the adjoint approach. In the latter case, the subdivision algorithm has the clear advantage over the reference point method since it depends on the accuracy of the state equation only.

\begin{table}[htbp]
\caption{Comparison of CPU time of different methods (Subdivision algorithm in parallel: CPU time = Number of CPUs times wall clock time).}
\begin{center}\footnotesize
\renewcommand{\arraystretch}{1.3}
\begin{tabular}{|c|c|c|c|c|}\hline
~ & \# boxes / & \# Function & \# Adjoint & CPU time [h] \\
	~ & points & evaluations & evaluations & ~ \\
	\hline
	Subdivision \eqref{eq:MOP-2D} & $1383$ & $\approx 8.3 \cdot 10^6$ & $0$ & $\approx 1132$ \\
	Subdivision Spline $\mathbb{R}^5$ & $6057$ & $\approx 11.1 \cdot 10^6$ & $0$ & $\approx 1512$ \\
	Subdivision Spline $\mathbb{R}^{10}$ & $1880$ & $\approx 42 \cdot 10^6$ & $0$ & $\approx 5875$ \\
	Reference point \eqref{eq:State2_a} -- \eqref{eq:Optimality2} & $256$ & $26781$ & $1002$ & $\approx 12.2$ \\
	\hline
\end{tabular}
\end{center}
\end{table}
\label{tab:CPU_time}
From a computational point of view, the reference point method is significantly faster, cf.~Table~\ref{tab:CPU_time}. This is not surprising since, in order to have a good representation of a box, a large number of sample points needs to be evaluated. This is already numerically expensive in 3D and more so for a spline based approximation, i.e.~we need to pay a price for global optimality and the derivative free approach. For future work, it is therefore advisable to utilize the gradient based version of the subdivision algorithm in order to decrease the computational cost, provided that high accuracy can be achieved for the gradient.

%% file: 06_Conclusion.tex
In this article, we present different approaches for numerically solving multiobjective optimal control problems involving boundary control of the Navier-Stokes equations. 
In order to reduce the computational effort, model order reduction via proper orthogonal decomposition and Galerkin projection is used to compute a surrogate reduced model of ODEs. Different multiobjective optimal control algorithms are introduced and their respective advantages and disadvantages are discussed. The subdivision algorithm yields a box covering of globally optimal Pareto sets and can also deal with problems where the set is disconnected. Since its applicability depends critically on both the decision space dimension and the numerical effort of function evaluations, 
it is restricted to moderate decision space dimensions. To solve optimal control problems, the control function therefore has to be represented by a low number of parameters, e.~g.~by a harmonic function or spline coefficients. The reference point method converts the multiobjective optimal control problem into a sequence of scalar optimal control problems that can be solved using well known procedures from single objective optimal control theory, also for high decision space dimensions. To incorporate gradient information, the optimality system is derived for the scalar problem. In this case, the accuracy of the approximated gradient is of great importance, which becomes evident when comparing the two optimality systems that were derived. Provided the gradient accuracy is sufficient, the Pareto set can be approximated very efficiently, using the solution from previous scalar problems as initial guesses. The method is of local nature, however, such that these initial guesses as well as the computation of the first Pareto point are important. Additionally, the selection of target points becomes tedious for higher objective space dimensions. 

The computed Pareto set yields considerably more information about the system than the solution of a single objective optimization. A decision maker can use this information to make well-founded decisions on how to control the system or even to devise adaptive control strategies reacting to changing priorities such as the need to save energy. 
There are numerous applications where multiobjective optimal control for systems described by PDEs is of great interest, e.g.~minimizing the drag and maximizing the lift of airplanes, minimizing the weight and maximizing the stability of structures, optimal mixing with minimal cost, control of flow patterns and vortices in HVAC applications, to name a few. Moreover, applications with higher Reynolds numbers are certainly of great interest \cite{BN15}, which adds difficulties especially concerning the model order reduction. A question that needs to be addressed in the future is the error control for the reduced model and the respective gradient in order to guarantee optimality for the PDE based problem. Several authors have addressed this for scalar optimization with reduced order models \cite{Fah00, KV02, VP05, BC08, GV13, Las14} and additional questions arise in the context of multiobjective optimal control.

%% file: Appendix.tex
\appendix
\section{Reduced Order Model}
\label{appendix:ROM}
Here we state the model coefficients for the reduced order model \eqref{eq:eROM_state}, see e.g.~\cite{Fah00, BCB05} for details. We refer to the time average of the modified flow field by $\boldsymbol{U}_m = \left\langle \widetilde{\boldsymbol{U}}(x, t) \right\rangle$.
\begin{align*}
\mathcal{A}_i = &-\left((\boldsymbol{U}_m \cdot \nabla ) \boldsymbol{U}_m, \boldsymbol{\psi}_i \right)_{L^2} - \frac{1}{Re} \left(\nabla \boldsymbol{U}_m, \nabla \boldsymbol{\psi}_i \right)_{L^2}, \\
\mathcal{B}_{ij} = &-\left((\boldsymbol{U}_m \cdot \nabla ) \boldsymbol{\psi}_j, \boldsymbol{\psi}_i \right)_{L^2} -\left((\boldsymbol{\psi}_j \cdot \nabla ) \boldsymbol{U}_m, \boldsymbol{\psi}_i \right)_{L^2} - \frac{1}{Re} \left(\nabla \boldsymbol{\psi}_i, \nabla \boldsymbol{\psi}_j \right)_{L^2},\\
\mathcal{Q}_{jik} = &-\left((\boldsymbol{\psi}_i \cdot \nabla ) \boldsymbol{\psi}_k, \boldsymbol{\psi}_j \right)_{L^2}, \\
\mathcal{D}_i = &-\left(\boldsymbol{U}_c, \boldsymbol{\psi}_i \right)_{L^2}, \\
\mathcal{E}_i = &-\left((\boldsymbol{U}_m \cdot \nabla ) \boldsymbol{U}_c, \boldsymbol{\psi}_i \right)_{L^2} -\left((\boldsymbol{U}_c \cdot \nabla ) \boldsymbol{U}_m, \boldsymbol{\psi}_i \right)_{L^2} - \frac{1}{Re} \left(\nabla \boldsymbol{U}_c, \nabla \boldsymbol{\psi}_i \right)_{L^2} , \\
\mathcal{F}_{ij} = &-\left((\boldsymbol{U}_c \cdot \nabla ) \boldsymbol{\psi}_j, \boldsymbol{\psi}_i \right)_{L^2} -\left((\boldsymbol{\psi}_j \cdot \nabla ) \boldsymbol{U}_c, \boldsymbol{\psi}_i \right)_{L^2}, \\
\mathcal{G}_i = &-\left((\boldsymbol{U}_c \cdot \nabla ) \boldsymbol{U}_c, \boldsymbol{\psi}_i \right)_{L^2}.
\end{align*}

\section{Derivation of the optimality system}
\label{appendix:OS}
Here we derive the optimality system \eqref{eq:State1} -- \eqref{eq:Optimality1} for the reference point method with the scalarized cost functional \eqref{eq:J_scalar}. In order to satisfy the necessary condition for optimality, we require that all variations of the Lagrange functional \eqref{eq:Lagrange1} are zero:
\begin{gather*}
	\delta L = \frac{\partial L}{\partial \alpha} \delta \alpha + \frac{\partial L}{\partial \dot\alpha} \delta \dot\alpha + \frac{\partial L}{\partial \lambda} \delta \lambda + \frac{\partial L}{\partial \gamma} \delta \gamma + \frac{\partial L}{\partial \dot\gamma} \delta \dot\gamma = 0 \\
	\Leftrightarrow \int_{t_0}^{t_e} \left( \frac{\partial \overline{J}}{\partial \alpha} + \lambda^\top \left( \mathcal{B} + \frac{\partial \mathcal{C}(\alpha)}{\partial \alpha} + \mathcal{F} \gamma \right) \right) \delta \alpha  + \left( \frac{\partial \overline{J}}{\partial \gamma} + \lambda^\top \left( \mathcal{E} + \mathcal{F}\alpha + 2 \mathcal{G} \gamma \right) \right) \delta \gamma + \\
	+ \left( \dot{\alpha} - \mathcal{A} - \mathcal{B} \alpha - \mathcal{C}(j\alpha) - \mathcal{D} \dot{\gamma} - (\mathcal{E} - \mathcal{F} \alpha) \gamma - \mathcal{G} \gamma^2 \right) \delta \lambda + \lambda^\top \mathcal{D} \delta \dot \gamma - \lambda^\top \delta \dot \alpha \ dt = 0.
\end{gather*}
Using partial integration, this leads to
\begin{gather*}
	\int_{t_0}^{t_e} \left( \dot \lambda^\top + \frac{\partial \overline{J}}{\partial \alpha} + \lambda^\top \left( \mathcal{B} + \frac{\partial \mathcal{C}(\alpha)}{\partial \alpha} + \mathcal{F} \gamma \right) \right) \delta \alpha + \\
	+ \left( \dot{\alpha} - \mathcal{A} - \mathcal{B} \alpha - \mathcal{C}(j\alpha) - \mathcal{D} \dot{\gamma} - (\mathcal{E} - \mathcal{F} \alpha) \gamma - \mathcal{G} \gamma^2 \right) \delta \lambda + \\
	+ \left( -\dot\lambda^\top \mathcal{D} + \frac{\partial \overline{J}}{\partial \gamma} + \lambda^\top \left( \mathcal{E} + \mathcal{F}\alpha + 2 \mathcal{G} \gamma \right) \right) \delta \gamma - \lambda^\top\delta\alpha \big|_{t_0}^{t_e} + \lambda^\top \mathcal{D} \delta \gamma \big|_{t_0}^{t_e} = 0.
\end{gather*}
Since we require each variation to be zero individually, we obtain the equations \eqref{eq:State1} -- \eqref{eq:Optimality1} and two boundary conditions at $t_0$ and $t_e$, respectively. We have a fixed value $\alpha(t_0) = \alpha_0$ such that $\delta \alpha (t_0) = 0$ and hence, we have no initial condition $\lambda(t_0)$. In contrast, $\alpha(t_e)$ is arbitrary such that $\delta \alpha (t_e) \neq 0$ which results in the terminal condition $\lambda(t_e) = 0$. Since we want to accept arbitrary initial values $\gamma(t_0)$ and $\lambda(t_0)$ is also free, we obtain an additional condition $\lambda^\top(t_0) \mathcal{D} = 0$. The second optimality system \eqref{eq:State2_a} -- \eqref{eq:Optimality2} is derived analogously.